\documentclass[a4paper]{article}
\usepackage{multicol,color}
\usepackage{amsmath,latexsym,amsbsy,amssymb}
\usepackage{enumerate}
\usepackage{pifont,bbding}
\usepackage[latin1]{inputenc}
\usepackage{amsfonts}
\def\disp{\displaystyle}

\def\Limsup{\mathop{{\rm Lim}\,{\rm sup}}}

\def\tto{\;{\lower 1pt\hbox{$\rightarrow$}}\kern -10pt
\hbox{\raise 2pt \hbox{$\rightarrow$}}\;}
\def\Hat{\widehat}
\def\Tilde{\widetilde}

\def\ra{\rangle}
\def\la{\langle}
\def\ve{\varepsilon}
\def\B{\mathbb{B}}
\def\h{\hfill\Box}
\def\R{\mathbb{R}}
\def\N{\mathbb{N}}

\def\ox{\bar{x}}
\def\oy{\bar{y}}
\def\oz{\bar{z}}
\def\ov{\bar{v}}
\def\ou{\bar{u}}

\def\ow{\bar{w}}

\def\ob{\bar{b}}
\def\t{\tau}

\def\gph{\mbox{\rm gph}\,}
\def\epi{\mbox{\rm epi}\,}
\def\dim{\mbox{\rm dim}\,}
\def\dom{\mbox{\rm dom}\,}

\def\h{\hfill\triangle}

\def\O{\Omega}

\def\ph{\varphi}
\def\emp{\emptyset}

\def\oR{\overline{\R}}

\def\lm{\lambda}
\def\gg{\gamma}
\def\dd{\delta}
\def\al{\alpha}
\def\bb{\beta}
\def\Th{\Theta}
\def\N{I\!\!N}
\def\th{\theta}
\def\vt{\vartheta}
\def\Hat{\widehat}

\def\Tilde{\widetilde}
\def\tilde{\widetilde}

\def\ra{\rangle}
\def\la{\langle}
\def\ve{\varepsilon}
\def\B{\mathbb{B}}
\def\h{\hfill\Box}
\def\R{\mathbb{R}}
\def\N{\mathbb{N}}

\def\t{\tau}
\def\ox{\bar{x}}
\def\op{\bar{p}}
\def\oy{\bar{y}}
\def\oz{\bar{z}}
\def\ov{\bar{v}}
\def\ou{\bar{u}}
\def\ow{\bar{w}}

\def\ob{\bar{b}}

\def\gph{\mbox{\rm gph}\,}
\def\epi{\mbox{\rm epi}\,}
\def\dim{\mbox{\rm dim}\,}
\def\dom{\mbox{\rm dom}\,}

\def\h{\hfill\triangle}

\def\O{\Omega}

\def\ph{\varphi}
\def\emp{\emptyset}

\def\oR{\overline{\R}}

\def\lm{\lambda}
\def\gg{\gamma}
\def\dd{\delta}
\def\al{\alpha}
\def\bb{\beta}
\def\Th{\Theta}
\newcounter{lk}

\title{Extended Euler-Lagrange and Hamiltonian Conditions in Optimal Control of Sweeping Processes with Controlled Moving Sets}
\author{Nguyen D. Hoang\footnote{Departamento de Ingenier\'ia Matem\'atica, Universidad de Concepci\'on, Casilla 160-C, Concepci\'on, Chile (hnguyen@ing-mat.udec.cl).}$\;\mbox{ and }\;$Boris S. Mordukhovich\footnote{Department of Mathematics, Wayne State University, Detroit, Michigan, USA (boris@math.wayne.edu). Research of this author was partly supported by the USA National Science Foundation under grant DMS-1512846 and by the USA Air Force Office of Scientific Research under grant \#15RT0462.}}
\begin{document}
\maketitle\vspace*{-0.2in}

\small{\bf Abstract.} This paper concerns optimal control problems for a class of sweeping processes governed by discontinuous unbounded differential inclusions that are described via normal cone mappings to controlled moving sets. Largely motivated by applications to hysteresis, we consider a general setting where moving sets are given as inverse images of closed subsets of finite-dimensional spaces under nonlinear differentiable mappings dependent on both state and control variables. Developing the method of discrete approximations and employing generalized differential tools of first-order and second-order variational analysis allow us to derive nondegenerated necessary optimality conditions for such problems in extended Euler-Lagrange and Hamiltonian forms involving the Hamiltonian maximization. The latter conditions of the Pontryagin Maximum Principle type are the first in the literature for optimal control of sweeping processes with control-dependent moving sets.


\newtheorem{Theorem}{Theorem}[section]
\newtheorem{Proposition}[Theorem]{Proposition}
\newtheorem{Remark}[Theorem]{Remark}
\newtheorem{Lemma}[Theorem]{Lemma}
\newtheorem{Corollary}[Theorem]{Corollary}
\newtheorem{Definition}[Theorem]{Definition}
\newtheorem{Example}[Theorem]{Example}
\newtheorem{Assumptions}[Theorem]{Assumptions}
\renewcommand{\theequation}{\thesection.\arabic{equation}}
\normalsize
\def\proof{\normalfont\medskip
{\noindent\itshape Proof.\hspace*{6pt}\ignorespaces}}
\def\endproof{$\h$\vspace*{0.1in}}\vspace*{-0.15in}

\section{Introduction}
\setcounter{equation}{0}

The basic sweeping process (``processus du rafle") was introduced by Moreau \cite{mor_frict} in the form
\begin{eqnarray}\label{sp}
\dot{x}(t)\in-N\big(x(t);C(t)\big)\;\mbox{ a.e. }\;t\in[0,T],
\end{eqnarray}
where $N(x;\O)$ stands for the normal cone to a convex set $\O\subset\R^n$ at $x$ defined by
\begin{eqnarray}\label{nor-conv}
N(x;\O):=\left\{\begin{array}{ll}
\big\{v\in\R^n\big|\;\la v,u-x\ra\le 0\;\mbox{ for all }\;u\in\O\big\}&\mbox{if }\;x\in\O,\\
\emp&\mbox{otherwise},
\end{array}\right.
\end{eqnarray}
and where the convex variable set $C(t)$ continuously evolves in time. It has been realized that the Cauchy problem $x(0)=x_0\in C(0)$ for \eqref{sp} admits a unique solution (see, e.g., \cite{ct}), and hence there is no sense to consider optimization problems for the sweeping differential inclusion \eqref{sp}. This is totally different from the developed optimal control theory for Lipschitzian differential inclusions of the type
\begin{eqnarray}\label{di}
\dot x(t)\in F\big(x(t)\big)\;\mbox{ a.e. }\;t\in[0,T],
\end{eqnarray}
which arose from the classical one for controlled differential equations
\begin{eqnarray}\label{de}
\dot x(t)=f\big(x(t),u(t)\big),\;\;u(t)\in U\;\mbox{ a.e. }\;t\in[0,T]
\end{eqnarray}
with $F(x):=f(x,U)=\{y\in\R^n|\;y=f(x,u)\;\mbox{ for some }\;u\in U\}$ in \eqref{di}; see, e.g., the books \cite{mordukhovich,vinter} with the references therein as well as more recent publications devoted to optimal control of \eqref{di}.

It was suggested in \cite{chhm}, probably for the first time in the literature, to formulate optimal control problems for \eqref{sp} by entering control functions into the moving sets $C(t)$ in \eqref{sp}, i.e., considering the moving set  {\em control parametrization} in the form
\begin{eqnarray}\label{cp}
C(t)=C\big(u(t)\big)\;\mbox{ for all }\;t\in[0,T]
\end{eqnarray}
with respect to some collections of admissible controls $u(\cdot)$ satisfying appropriate constraints. In this way we arrive at new and very challenging classes of optimal control problems on minimizing certain Bolza-type cost functionals over feasible solutions to {\em highly non-Lipschitzian} unbounded differential inclusions under {\em irregular pointwise} state-control constraints
\begin{eqnarray}\label{state-con}
x(t)\in C\big(u(t)\big)\;\mbox{ for all }\;t\in[0,T],
\end{eqnarray}
which intrinsically arise from \eqref{sp} and \eqref{cp} due the normal cone construction in \eqref{nor-conv}. It occurs that not only results but also methods developed in optimal control theory for controlled differential equations \eqref{de} and Lipschitzian differential inclusions \eqref{di} are not suitable for applications to the new classes of sweeping control systems that appear in this way. Papers \cite{chhm,h1} present significant extensions to sweeping control systems of type \eqref{sp}, \eqref{cp}, and \eqref{state-con} of the {\em method of discrete approximations} developed in \cite{m95,mordukhovich} for Lipschitzian differential inclusions of type \eqref{di}. Major new ideas in the obtained extensions consist of marring the discrete approximation approach to recently established {\em second-order subdifferential calculus} and {\em explicit computations} of the corresponding second-order constructions of variational analysis. The strongest results established by such a device in \cite{h1} concern necessary optimality conditions for the generalized Bolza problem with the controlled sweeping dynamics in \eqref{sp}, \eqref{cp}, and \eqref{state-con} described by the moving convex polyhedra of the type
\begin{eqnarray}\label{sw-con2}
C(t):=\big\{x\in\R^n\big|\;\la u_i(t),x\ra\le b_i(t),\;i=1,\ldots,m\big\},\quad t\in[0,T],
\end{eqnarray}
where both actions $u_i(t)$ and $b_i(t)$ are involved in control. Other developments of the discrete approximation approach to derive necessary conditions for controlled sweeping systems with controls not only in the moving sets \eqref{cp} but also in additive perturbations of \eqref{sp} are given in \cite{cm1}--\cite{cm3}, where the reader can find applications of the obtained results to the practical crowd motion model of traffic equilibrium. The method of discrete approximations  was also implemented in \cite{dfm} to study various optimal control issues for evolution inclusions governed by one-sided Lipschitzian mappings and in \cite{h2} for those described by maximal monotone operators in Hilbert spaces, but without deriving necessary optimality conditions.

Note that the necessary optimality conditions obtained in the previous papers \cite{cm1}--\cite{h1} do not contain the formalism of the Pontryagin Maximum Principle (PMP) \cite{pontryagin} (i.e., the maximization of the corresponding Hamiltonian function) established in classical optimal control of \eqref{de} and then extended to optimal control problems for Lipschitzian differential inclusions of type \eqref{di}.

To the best of our knowledge, necessary optimality conditions involving the maximization of the corresponding Hamiltonian were first obtained for sweeping control systems in \cite{bk}, where the authors considered a sweeping process with a strictly smooth, convex, and solid set $C(t)\equiv C$ in \eqref{sp} while with control functions entering linearly an adjacent ordinary differential equation. Further results with the maximum condition for global (as in \cite{bk}) minimizers were derived in \cite{ac} for the sweeping control system
\begin{eqnarray}\label{sp1}
\dot x(t)\in f\big(x(t),u(t)\big)-N\big(x(t);C(t)\big)\;\mbox{ a.e. }\;t\in[0,T],
\end{eqnarray}
where measurable controls $u(t)$ enter the additive smooth term $f$ while the uncontrolled moving set $C(t)$ is compact, uniformly prox-regular regular (close enough to convexity), and possesses a ${\cal C}^3$-smooth boundary for each $t\in[0,T]$ under some other assumptions. The very recent paper \cite{pfs} also concerns a (generally nonautonomous) sweeping control system in form \eqref{sp1} and derives necessary optimality conditions of the PMP type for global minimizers provided that the convex, solid, and compact set $C(t)\equiv C$ therein is defined by $C:=\{x\in\R^n|\;\psi(x)\le 0\}$ via a ${\cal C}^2$-smooth function $\psi$ under other assumptions, which are partly differ from \cite{ac}. The {\em penalty-type} approximation methods developed in \cite{ac}, \cite{bk}, and \cite{pfs} are different from each other, significantly based on the {\em smoothness} of uncontrolled moving sets while being totally distinct from the method of discrete approximations employed in our previous papers and in what follows.\vspace*{0.02in}

This paper addresses sweeping control systems modeled as
\begin{eqnarray}\label{evo_equa}
\dot x(t)\in f\big(t,x(t)\big)-N\big(g(x(t));C(t,u(t))\big)\;\mbox{ a.e. }\;t\in[0,T],\quad x(0)=x_0\in C\big(0,u(0)\big),
\end{eqnarray}
where the controlled moving set is given by
\begin{eqnarray}\label{mov-set}
C(t,u):=\big\{x\in\R^n\big|\;\psi(t,x,u)\in\Th\big\},\quad(t,u)\in[0,T]\times\R^m,
\end{eqnarray}
with $f\colon[0,T]\times\R^n\to\R^n$, $g\colon\R^n\to\R^n$, $\psi\colon[0,T]\times\R^n\times\R^m\to\R^s$, and $\Th\subset\R^s$. Definition~\eqref{mov-set} amounts to saying that $C(t,u)$ is the {\em inverse image} of the set $\Th$ under the mapping $x\mapsto\psi(t,g(x),u)$ for any $((t,u)$. Throughout the paper we assume that the set $\Th$ is {\em locally closed} around the reference point. We do not impose any convexity of $C(t,u)$ and use in \eqref{evo_equa} the (basic, limiting, Mordukhovich) {\em normal cone} to an arbitrary locally closed set $\O\subset\R^n$ at $\ox\in\R^n$ defined by
\begin{eqnarray}\label{nc}
N(\ox;\O):=\left\{\begin{array}{ll}
\big\{v\in\R^n\big|\;\exists x_k\to\ox,\;\al_k\ge 0,\;w_k\in\Pi(x_k;\O),\;\al_k(x_k-w_k)\to v&\mbox{if }\;\ox\in\O,\\
\emp&\mbox{otherwise},
\end{array}\right.
\end{eqnarray}
where $\Pi(x;\O)$ stands for the Euclidean projector of $x$ onto $\O$. When $\O$ is convex, the normal cone \eqref{nc} reduces to  the one \eqref{nor-conv} in the sense of convex analysis, but in general the multifunction $x\tto N(x;\O)$ is nonconvex-valued while satisfying a {\em full calculus} together with the associated subdifferential of extended-real-valued functions and coderivative of set-valued mappings considered below. Such a calculus is due to {\em variational/extremal principles} of variational analysis; see \cite{mordukhovich,mord,rw} for more details.

Our major goal here is to study the optimal control problem $(P)$ of minimizing the cost functional
\begin{eqnarray}\label{eq:MP}
{\rm minimize}\;J[x,u]:=\varphi\big(x(T)\big)+\int_0^T\ell\big(t,x(t),u(t),\dot x(t),\dot u(t)\big)dt
\end{eqnarray}
over absolutely continuous control actions $u(\cdot)$ and the corresponding absolutely continuous trajectories $x(\cdot)$ of the sweeping differential inclusion \eqref{evo_equa} generated by the controlled moving set \eqref{mov-set}. It follows from \eqref{evo_equa} and the normal cone definition \eqref{nc} that the optimal control problem in \eqref{evo_equa} and \eqref{eq:MP} intrinsically contains the {\em pointwise constraints} on both {\em state} and {\em control} functions
\begin{eqnarray*}\label{const}
\psi\big(t,g(x(t)),u(t)\big)\in\Th\;\mbox{ for all }\;t\in[0,T].
\end{eqnarray*}
Note that the optimal control problem studied in \cite{h1} is a particular case of our problem $(P)$ that corresponds to the choice of $g(x):=x$ (identity operator), $\psi(t,x,u):=Ax-b$ and $Z:=\R^m_-$ in \eqref{evo_equa} and \eqref{mov-set}. Besides being attracted by challenging mathematical issues, our interest to more general sweeping control problems considering in this paper is largely motivated by applications to {\em rate-independent operators} that frequently appear, e.g., in various plasticity models and in the study of hysteresis. We discuss these and related topics in more details in Section~\ref{application-example} and also will devote a separate paper to such applications.

While the underlying approach to derive necessary optimality conditions for local minimizers of the above problem $(P)$ is the usage of the {\em method of discrete approximations} and {\em generalized differentiation}, similarly to \cite{h1} and our other publications on sweeping optimal control, some important elements of our technique here are significantly different from the previous developments. From one side, the new/modified technique allows us to establish {\em nondegenerated} necessary optimality conditions for local minimizers of $(P)$ in the {\em extended Euler-Lagrange} form for more general sweeping systems with relaxing several restrictive technical assumptions of \cite{h1} in its polyhedral setting. On the other hand, we obtain optimality conditions of the {\em Hamiltonian/PMP} type, which are new even for polyhedral moving sets as in \cite{h1} under an additional surjectivity assumption. In fact, the optimality conditions in the PMP form are the {\em first results} of this type for sweeping process with controlled moving sets.\vspace*{0.02in}

The rest of the paper is organized as follows. In Section~\ref{prelim} we formulate and discuss our standing assumptions and present necessity preliminaries from first-order and second-order generalized differentiation that are widely used for deriving the main results of the paper.

Section~\ref{discrete-approximations} concerns discrete approximations of feasible and local optimal solutions to the sweeping control problem $(P)$ with the verification of the required strong convergence. In Section~\ref{optimality-conditions} we derive the extended Euler-Lagrange conditions for local optimal solutions to $(P)$ by passing to the limit from discrete approximations and using the second-order subdifferential calculations.

Section~\ref{hamiltonian} contains necessary optimality conditions of the PMP type involving the maximization of the new Hamiltonian function, discusses relationships with the conventional Hamiltonians, and presents an example showing that the maximum principle in the conventional Hamiltonian form fails in our framework. More examples of some practical meaning in the areas of elastoplasticity and hysteresis are given in Section~\ref{application-example}. The final Section~\ref{conclusion} discusses some directions of the future research.

Throughout the paper we use standard notation of variational analysis and control theory; see, e.g., \cite{mordukhovich,rw,vinter}. Recall that $\N:=\{1,2,\ldots\}$, that $A^*$ stands for the transposed/adjoint matrix to $A$, and that $\B$ denotes the closed unit ball of the space in question.\vspace*{-0.2in}

\section{Standing Assumptions and Preliminaries}\label{prelim}
\setcounter{equation}{0}

Let us first formulate the major assumptions on the given data of problem $(P)$ that are standing throughout the whole paper. Since our approach to derive necessary optimality conditions for $(P)$ is based on the method of discrete approximations, we impose the a.e.\ continuity of the functions involved with respect to the time variable, although it is not needed for results dealing with discrete systems before passing to the limit. Note also that the time variable is never included in subdifferentiation. As mentioned above, the constraint set $Z$ in \eqref{mov-set} is assumed to be locally closed unless otherwise stated.\vspace*{0.05in}

Our {\em standing assumptions} are as follows:\vspace*{0.05in}

{\bf (H1)} There exits $L_f>0$ such that $\|f(t,x)-f(t,y)\|\le L_f\|x-y\|$ for all $x,y\in\R^n,\;t\in[0,T]$ and the mapping $t\mapsto f(t,x)$ is a.e.\ continuous on $[0,T]$ for each $x\in\R^n$.

{\bf(H2)} There exits $L_g>0$ such that $\|g(x)-g(y)\|\le L_g\|x-y\|$ for all $x,y\in\R^n$.

{\bf(H3)} For each $(t,u)\in[0,T]\times\R^m$, the mapping $\psi_{t,u}(x):=\psi(t,x,u)$ is ${\cal C}^2$-smooth around the reference points with the surjective derivative $\nabla\psi_{t,u}(x)$ satisfying
\begin{eqnarray*}
\|\nabla\psi_{t,u}(x)-\nabla\psi_{t,v}(x)\|\le L_{\psi}\|u-v\|
\end{eqnarray*}
with the uniform Lipschitz constant $L_\psi$. Furthermore, the mapping $t\mapsto\psi(t,x)$ is a.e.\ continuous on $[0,T]$ for each $x\in\R^n$ and $u\in\R^m$.

{\bf(H4)} There are a number $\tau>0$ and a mapping $\vt\colon\R^n\times\R^n\times\R^n\times\R^m\to\R^m$ locally Lipschitz continuous and uniformly bounded on bounded sets such that for all $t\in[0,T]$, $\ov\in N(\psi_{(t,\ou)}(\ox);\Th)$, and $x\in\psi^{-1}_{(t,u)}(\Th)$ with $u:=\ou+\vt(x-\ox,x,\ox,\ou)$ there exists $v\in N(\psi_{(t,u)}(x);\Th)$ satisfying $\|v-\ov\|\le\tau\|x-\ox\|$.

{\bf(H5)} The cost functions $\ph\colon\R^n\to\oR:=[-\infty,\infty)$ and $\ell(t,\cdot)\colon\R^{2(n+m)}\to\oR$ in \eqref{eq:MP} are bounded from below and lower semicontinuous (l.s.c.) around a given feasible solution to $(P)$ for a.e.\ $t\in[0,T]$, while the integrand $\ell$ is a.e.\ continuous in $t$ and is uniformly majorized by a summable function on $[0,T]$.\vspace*{0.05in}

Assumption (H4) is technical and seems to be the most restrictive. Let us show nevertheless that it holds automatically in the polyhedral setting of \cite{h1} and also for nonconvex moving sets.\vspace*{-0.05in}

\begin{Proposition}{\bf(validity of (H4) for controlled polyhedra).}\label{poly} Let
\begin{eqnarray*}
\psi\big(t,x,(u,b)\big):=\la x,u\ra-b\;\mbox{ and }\;\Th=\R^m_-
\end{eqnarray*}
in \eqref{mov-set}. Then condition ${\rm(H4)}$ is satisfied.
\end{Proposition}\vspace*{-0.05in}
{\bf Proof.} Pick $\ov\in N(\la\ox,\ou\ra-\ob;\R^m_-)$, $x\in\R^n$ and denote $\vt(x,y,z,u):=(0,\langle x,u\rangle)$. Choose $(u,b):=(\ou,\ob)+\vt(x-\ox,x,\ox,\ou)$ and hence get $u=\ou$ and $b=\ob+\langle x-\ox,\ou\rangle$, which results in $\la\ox,\ou\ra-\ob=\la x,u\ra-b$. Then $N(\la\ox,\ou\ra-\ob;\R^m_-)=N(\la x,u\ra-b;\R^m_-)$. We can choose $v:=\ov$, and thus condition (H4) is satisfied with $v:=\ov$ for any number $\tau\ge 1$ therein. $\h$\vspace*{0.05in}

The following simple example illustrates that (H4) is also satisfied in standard nonconvex settings.\vspace*{-0.05in}

\begin{Example}{\bf(validity of (H4) for nonconvex moving sets).}\label{h4n} {\rm Consider the nonconvex set
\begin{eqnarray*}
C(t,u)=\big\{(x\in\R\big|\;x^2\ge-u+1\big\},
\end{eqnarray*}
which corresponds to $\psi(x,u):=x^2+u-1$ and $\Th:=[0,\infty)$ in \eqref{mov-set}. To verify (H4) in this setting, denote $\vt(x,y,z,u):=-x(y+z)$ and pick any
$\ov\in N(\ox^2+\ou-1;[0,\infty))$ and $x\in\R^n$. Choosing now $u:=\ou-(x-\ox)(x+\ox)$, we get $x^2+u-1=\ox^2+\ou-1$ and thus verify (H4) with $v:=\ov$ for every $\tau\ge 1$.}
\end{Example}\vspace*{-0.05in}

Let us next discuss condition (H3), which plays a significant role is deriving some major results of the paper. This condition, which is equivalent in the finite-dimensional setting under consideration to the full rank of the Jacobian matrix $\nabla\psi_{t,u}(x)$, amounts to {\em metric regularity} of the mapping $x\mapsto\psi_{t,u}(x)$ by the seminal Lyusternik-Graves theorem; see, e.g., \cite[Theorem~1.57]{mordukhovich}. The following normal cone calculus rule is a consequence of \cite[Theorem~1.17]{mordukhovich}.\vspace*{-0.05in}

\begin{Proposition}{\bf(normal cone representation for inverse images).}\label{nor_inve} Under the validity of {\rm(H3)} the normal cone \eqref{nc} to the controlled moving set \eqref{mov-set} is represented by
\begin{eqnarray*}
N\big(x;C(t,u)\big)=\nabla\psi_{t,u}(x)^*N\big(\psi_{t,u}(x);\Th\big)\;\mbox{ whenever }\;\psi(t,x,u)\in\Th.
\end{eqnarray*}
\end{Proposition}

To proceed further, we recall some constructions of first-order and second-order generalized differentiation for functions and multifunctions/set-valued mappings needed in what follows; see \cite{mordukhovich,mord} for detailed expositions. All these constructions are generated geometrically by our basic normal cone \eqref{nc}.

Given a set-valued mapping $F\colon\R^n\tto\R^q$ and a point $(\ox,\oy)\in\gph F$ from its graph
\begin{eqnarray*}
\gph F:=\big\{(x,y)\in\R^n\times\R^q\big|\;y\in F(x)\big)\big\},
\end{eqnarray*}
the {\em coderivative} $D^*F(\ox,\oy)\colon\R^q\tto\R^n$ of $F$ at $(\ox,\oy)$ is defined by
\begin{eqnarray}\label{cod}
D^*F(\ox,\oy)(u):=\big\{v\in\R^n\big|\;(v,-u)\in N\big((\ox,\oy);\gph F\big)\big\},\quad u\in\R^q,
\end{eqnarray}
where $\oy$ is omitted in the notation if $F\colon\R^n\to\R^q$ is single-valued. If furthermore $F$ is ${\cal C}^1$-smooth around $\ox$ (or merely strictly differentiable at this point), we have $D^*F(\ox)(v)=\{\nabla F(\ox)^*v\}$ via the adjoint Jacobian matrix. In general, the coderivative \eqref{cod} is a positively homogeneous multifunction satisfying comprehensive calculus rules and providing complete characterizations of major well-posedness properties in variational analysis related to Lipschitzian stability, metric regularity, and linear openness; see \cite{mordukhovich,rw}.

For an extended-real-valued function $\phi\colon\R^n\to\oR$ finite at $\ox$, i.e., with $\ox\in\dim\phi$, the (first-order) {\em subdifferential} of $\phi$ at $\ox$ is defined geometrically by
\begin{eqnarray}\label{1sub}
\partial\phi(\ox):=\{v\in\R^n\big|\;(v,-1)\in N\big((\ox,\phi(\ox));\epi\phi\big)\big\}
\end{eqnarray}
via the normal cone \eqref{nc} to the epigraphical set $\epi\phi:=\{(x,\al)\in\R^{n+1}|\;\al\ge\phi(x)\}$. If $\phi(x):=\dd_\O(x)$, the indicator function of a set $\O$ that equals to 0 for $x\in\O$ and to $\infty$ otherwise, we get $\partial\phi(\ox)=N(\ox;\O)$. Given further $\ov\in\partial\phi(\ox)$, the {\em second-order subdifferential} (or generalized Hessian) $\partial^2\phi(\ox,\ov)\colon\R^n\tto\R^n$ of $\phi$ at $\ox$ relative to $\ov$ is defined as the coderivative of the first-order subdifferential by
\begin{eqnarray}\label{2sub}
\partial^2\phi(\ox,\ov)(u):=(D^*\partial\phi)(\ox,\ov)(u),\quad u\in\R^n,
\end{eqnarray}
where $\ov=\nabla\phi(\ox)$ is omitted when $\phi$ is differentiable at $\ox$. If $\phi$ is ${\cal C}^2$-smooth around $\ox$, then \eqref{2sub} reduces to the classical (symmetric) Hessian matrix
\begin{eqnarray*}
\partial^2\phi(\ox)(u)=\big\{\nabla^2\phi(\ox)u\big\}\;\mbox{ for all }\;u\in\R^n.
\end{eqnarray*}

For applications in this paper we also need partial versions of the above subdifferential constructions for functions of two variables $\phi\colon\R^n\times\R^m\to\oR$. Consider the {\em partial first-order subdifferential} mapping $(x,w)\mapsto\partial_x\phi(x,w)$ for $\ph(x,w)$ with respect to $x$ by
\begin{eqnarray*}\label{1par}
\partial_x\phi(x,w):=\big\{\;\mbox{set of subgradients }\;v\in\R^n\;\mbox{ of }\;\phi_w:=\phi(\cdot,w)\;\mbox{ at }\;x\big\}=\partial\phi_w(x)
\end{eqnarray*}
and then, picking $(\ox,\ow)\in\dom\phi$ and $\ov\in\partial_x\phi(\ox,\ow)$, define the {\em partial second-order subdifferential} of $\phi$ with respect to $x$ at $(\ox,\ow)$ relative to $\ov$ by
\begin{eqnarray}\label{2par}
\partial^2_x\phi(\ox,\ow,\ov)(u):=\big(D^*\partial_x\phi)(\ox,\ow,\ov)(u)\;\mbox{ for all }\;u\in\R^n.
\end{eqnarray}
If $\phi$ is ${\cal C}^2$-smooth around $(\ox,\ow)$, we have the representation
\begin{eqnarray*}
\partial^2\phi(\ox,\ow)(u)=\big\{\big(\nabla^2_{xx}\phi(\ox,\ow)^*u,\nabla^2_{xw}\phi(\ox,\ow)^*u\big)\big\}\;\mbox{ for all }\;u\in\R^n.
\end{eqnarray*}

Taking into account the controlled moving set structure \eqref{mov-set}, important roles in this paper are played by the parametric constraint system
\begin{eqnarray}\label{ctild}
S(w):=\big\{x\in\R^n\big|\;\psi(x,w)\in\Th\big\},\quad w\in\R^m,
\end{eqnarray}
and the {\em normal cone mapping} ${\cal N}\colon\R^n\times\R^m\tto\R^n$ associated with \eqref{ctild} by
\begin{eqnarray}\label{nm}
{\cal N}(x,w):=N\big(x;S(w)\big)\;\mbox{ for }\;x\in S(w).
\end{eqnarray}
It is easy to see that the mapping ${\cal N}$ in \eqref{nm} admits the composite representation
\begin{eqnarray}\label{2comp}
{\cal N}(x,w)=\partial_x\phi(x,w)\;\mbox{ with }\;\phi(x,w):=\big(\dd_{\Th}\circ\psi\big)(x,w)
\end{eqnarray}
via the ${\cal C}^2$-smooth mapping $\psi\colon\R^n\times\R^m\to\R^s$ from \eqref{ctild} and the indicator function $\dd_\Th$ of the closed set $\Th\subset\R^s$. It follows directly from \eqref{2comp} due to the second-order subdifferential construction \eqref{2par} that
\begin{eqnarray}\label{nm1}
\partial^2_x\phi(\ox,\ow,\ov)(u)=D^*{\cal N}(\ox,\ow,\ov)(u)\;\mbox{ for any }\;\ov\in{\cal N}(\ox,\ow)\;\mbox{ and }\;u\in\R^n.
\end{eqnarray}
Applying now the second-order chain rule from \cite[Theorem~3.1]{BR} to the composition in \eqref{nm1} allows us to compute the coderivative of the normal cone mapping \eqref{nm} via the given data of \eqref{ctild}.\vspace*{-0.06in}

\begin{Proposition}{\bf(coderivative of the normal cone mapping for inverse images).}\label{morout} Assume that $\psi$ is ${\cal C}^2$-smooth around $(\ox,\ow)$, and that the partial Jacobian matrix $\nabla_x\psi(\ox,\ow)$ is of full rank. Then for each $\ov\in{\cal N}(\ox,\ow)$ there is a unique vector $\op\in N_{\Th}(\psi(\bar{x},\bar{w})):=N(\psi(\bar{x},\bar{w});\Th)$ satisfying
\begin{eqnarray}\label{unique_re}
\nabla_{x}\psi(\bar{x},\bar{w})^*\op=\bar{v}
\end{eqnarray}
and such that the coderivative of the normal cone mapping is computed for all $u\in\R^n$ by
\begin{eqnarray*}
D^*{\cal N}(\bar{x},\bar{w},\bar{v})(u)=\left[\begin{array}{c}
\nabla_{xx}^{2}\langle\bar p,\psi\rangle(\bar{x},\bar{w})\\\nabla_{xw}^{2}\langle\bar p,\psi\rangle(\bar{x},\bar{w})
\end{array}\right]u+\nabla\psi(\bar{x},\bar{w})^*D^*N_{\Th}\big(\psi(\bar{x},\bar{w}),\bar p\big)\big(\nabla_{x}\psi(\bar{x},\bar{w})u\big).
\end{eqnarray*}
\end{Proposition}\vspace*{-0.02in}

Thus Proposition~\ref{morout} reduces the computation of $D^*{\cal N}$ to that of $D^*N_{\Th}$, which has been computed via the given data for broad classes of sets $\Th$; see, e.g., \cite{heoutsur,mord,MO07,BR} for more details and references.\vspace*{-0.15in}

\section{Discrete Approximations}\label{discrete-approximations}
\setcounter{equation}{0}

In this section we construct a well-posed sequence of discrete approximations for feasible solutions to the constrained sweeping dynamics in \eqref{evo_equa}, \eqref{mov-set} and for local optimal solutions to the sweeping optimal control problem $(P)$. The results obtained here establish the $W^{1,2}$-strong convergence of discrete approximations while being free of generalized differentiation. They are certainly of independent interest from just their subsequent applications to deriving necessary optimality conditions for problem $(P)$.\vspace*{0.02in}

Starting with discrete approximations of the sweeping differential inclusion \eqref{evo_equa}, we replace the time derivative therein by the Euler finite difference $\dot x(t)\approx[x(t+h)-x(t)]/h$ and proceed as follows. For each $k\in\N$ define $h_k:=T/k$ and consider the discrete mesh $T_k:=\{t_j^k:=jh_k|\,j=0,1,\ldots,k\}$. Then the sequence of discrete approximations of \eqref{evo_equa} is given by
\begin{eqnarray}\label{evo_equa_discrete}
\frac{x^k_{j+1}-x^k_j}{h_k}\in f(t^k_j,x^k_{j})-N\big(g(x_j^k);C(t^k_j,u_j^k)\big),\;j=0,\ldots,k-1;\;x^k_0=x_0\in C\big(0,u(0)\big).
\end{eqnarray}

The following major result on a {\em strong approximation} of {\em feasible solutions} to the controlled sweeping process in \eqref{evo_equa}, \eqref{mov-set} by feasible solutions to the discretized systems in \eqref{evo_equa_discrete} is essentially different from the related one in \cite[Theorem~3.1]{h1}. Besides being applied to more general systems, it eliminates or significantly relax some restrictive technical assumptions imposed in \cite{h1} for the polyhedral controlled sets \eqref{sw-con2} by using another proof technique under the surjectivity assumption in (H3). Note furthermore that the choices of the reference feasible pair $(\ox(\cdot),\ou(\cdot))$ in \cite{h1} and Theorem~\ref{feasible-approx} below are also different: instead of the actual choice of $(\ox(\cdot),\ou(\cdot))\in W^{2,\infty}[0,T]\times W^{2,\infty}[0,T]$ we now have the less restrictive pick $(\ox(\cdot),\ou(\cdot))\in{\cal C}^1[0,T]\times{\cal C}[0,T]$ and establish its strong approximation in the $W^{1,2}\times{\cal C}$ topology instead of $W^{1,2}\times W^{1,2}$ in \cite{h1}. This actually effects the types of local minimizers for which we derive necessary optimality conditions in the setting of this paper; see Definition~\ref{ilm} below. If however the mapping $g$ is linear in \eqref{evo_equa} and if $\ou(\cdot)\in W^{1,2}[0,T]$, we obtain the same $W^{1,2}$-approximation of the reference control as in the polyhedral framework of \cite{h1}.\vspace*{-0.05in}

\begin{Theorem}{\bf(strong discrete approximation of feasible sweeping solutions).}\label{feasible-approx} Let the pair $(\ox(\cdot),\ou(\cdot))\in{\cal C}^1[0,T]\times{\cal C}[0,T]$ satisfy \eqref{evo_equa} with the moving set $C(t,u)$ from \eqref{mov-set} under the validity of the standing assumptions in {\rm(H1)--(H4)}. Then there exists a sequence $\{(x^k(\cdot),u^k(\cdot))\}$ of piecewise linear functions on $[0,T]$ satisfying the discrete inclusions \eqref{evo_equa_discrete} with $(x^k(0),u^k(0))=(\ox_0,\ou(0))$ and such that $\{(x^k(\cdot),u^k(\cdot))\}$ converges to $(\ox(\cdot),\ou(\cdot))$ in the norm topology of $W^{1,2}([0,T];\R^{n})\times{\cal C}([0,T];\R^m)$. If in addition the mapping $g(\cdot)$ in \eqref{evo_equa} is linear and $\ou(\cdot)\in W^{1,2}([0,T];\R^m)$, then the sequence $\{(x^k(\cdot),u^k(\cdot))\}$ converges to $(\ox(\cdot),\ou(\cdot))$ in the norm topology of $W^{1,2}([0,T];\R^n)\times W^{1,2}([0,T];\R^m)$.
\end{Theorem}\vspace*{-0.05in}
{\bf Proof.} Fix $k\in\N$, choose $x_0^k:=x_0$ and $u^k_0:=\ou(0)$, and then construct $(x_j^k,u_j^k)$ for $j=1,\ldots,k$ by induction. Suppose that $x_j^k$ is known and satisfies $\|x_j^k-\ox(t_j^k)\|\le 1$ without loss of generality. Define
\begin{eqnarray}\label{ep}
\epsilon_j:=\|\ox(t^k_{j+1})-\ox(t^k_j)-h_k\dot x(t^k_j)\|\;\mbox{ for all }\;j=0,\ldots,k-1
\end{eqnarray}
and deduce from the validity of \eqref{evo_equa} at $t_j^k$ that $-\dot\ox(t_j^k)+f(t_j^k,\ox(t_j^k))\in N(g(\ox(t_j^k);C(t_j^k,\ou(t_j^k)))$. Using
$\ox(\cdot)\in{\cal C}^1([0,T];\R^n)$ allows us to find $\eta>0$ such that
\begin{eqnarray*}
\|\nabla\psi_{t_j^k,\ou(t_j^k)}\|\le\eta\;\mbox{ and }\;\|-\dot\ox(t_j^k)+f\big(t_j^k,\ox(t_j^k)\big)\|\le\eta\;\mbox{ whenever }\;k\in\N\;\mbox{ and }\;j=0,\ldots,k.
\end{eqnarray*}
The surjectivity of $\nabla\psi_{t_j^k,\ou(t_j^k)}$ ensures by the open mapping theorem the existence of $M>0$ for which
\begin{eqnarray*}
\B\subset\big(\nabla\psi_{t_j^k,\ou(t_j^k)}\big)^*(M\B).
\end{eqnarray*}
Combining it with Proposition~\ref{nor_inve} tells us that
\begin{eqnarray*}
N\big(g(\ox(t_j^k);C(t_j^k,\ou(t_j^k))\big)\cap\eta\B=\big(\nabla\psi_{t_j^k,\ou(t_j^k)}\big)^*\Big(N\big(\psi_{t_j^k,\ou(t_j^k)}(g(\ox(t_j^k));\Th\big)\cap\eta M\B\Big).
\end{eqnarray*}
Since $-\dot\ox(t_j^k)+f(t_j^k,\ox(t_j^k))\in N(g(\ox(t_j^k);C(t_j^k,\ou(t_j^k)))\cap\eta\B$, we find $w\in N(\psi_{t_j^k,\ou(t_j^k)}(g(\ox(t_j^k)));\Th)$ with $\|w\|\le\eta M$ satisfying the equality
\begin{eqnarray*}
-\dot\ox(t_j^k)+f\big(t_j^k,\ox(t_j^k)\big)=\nabla\psi_{t_j^k,\ou(t_j^k)}^*w.
\end{eqnarray*}
Using now the mapping $\vt(\cdot)$ from (H4) gives us vectors $u_j^k$ and $\tilde w\in N(\psi_{t_j^k,\ou_j^k}(x_j^k);\Th)$ for which
\begin{eqnarray*}\label{d}
u_j^k=\ou(t_j^k)+d\big(g(x_j^k)-g(\ox(t_j^k)),g(x_j^k),g(\ox(t_j^k)),\ou(t_j^k)\big)\;\mbox{ and }\;\|w-\tilde w\|\le\tau\|g(x_j^k)-g(\ox(t_j^k))\|.
\end{eqnarray*}
By the assumed uniform boundedness of the mapping $\vt(\cdot)$ it is easy to adjust $\tau>0$ so that $\|u_j^k-\ou(t_j^k)\|\le\tau\|g(x_j^k)-g(\ox(t_j^k))\|$. Denoting $v_j^k:=\big(\nabla\psi_{t_j^k,u_j^k}\big)^*\tilde w$, we get $v_j^k\in N(g(x_j^k);C(t_j^k,u_j^k))$ by employing Proposition~\ref{nor_inve} and then arrive at the following estimates:
\begin{align*}
\big\|v_j^k-\big(-\dot\ox(t_j^k)+f(t_j^k,\ox(t_j^k))\big)\big\|&=\big\|\big(\nabla\psi_{t_j^k,\ou(t_j^k)}\big)^*w-\big(\nabla\psi_{t_j^k,u_j^k}\big)^*\tilde w\big\|\\&\le\big\|\big(\nabla\psi_{t_j^k,\ou(t_j^k)}-\nabla\psi_{t_j^k,u_j^k}\big)^*w\big\|+\big\|\big(\nabla\psi_{t_j^k,u_j^k}\big)^*(w-\tilde w)\big\|\\
&\le\eta M L_{\psi}\|u_j^k-\ou(t_j^k)\|+\tau L_g\|x_j^k-\ox(t_j^k)\|\big(\eta+L_{\psi}\|u_j^k-\ou(t_j^k)\|\big)\\
&\le\Big(\eta M L_{\psi}\tau L_g+\tau L_g\big(\eta+L_{\psi}\tau L_g\|x_j^k-\ox(t_j^k)\|\big)\Big)\|x_j^k-\ox(t_j^k)\|\\
&\le\Big(\eta M L_{\psi}\tau L_g+\tau L_g\eta+L_{\psi}\tau^2L^2_g\|x_j^k-\ox(t_j^k)\|\Big)\|x_j^k-\ox(t_j^k)\|.
\end{align*}
Denoting further $\alpha:=\eta M L_{\psi}\tau L_g+\tau\eta L_g$ and $\beta:=L_{\psi}\tau^2 L^2_g$, we get from the above that
\begin{eqnarray}\label{derivative_est}
\big\|v_j^k-\big(-\dot\ox(t_j^k)+f(t_j^k,\ox(t_j^k))\big)\big\|\le\big(\alpha+\beta\|x_j^k-\ox(t_j^k)\|\big)\|x_j^k-\ox(t_j^k)\|.
\end{eqnarray}
Now we are ready to construct the next iterate $x^k_{j+1}$ by
\begin{eqnarray*}
x^k_{j+1}:=x^k_j+h_k f(t^k_j,x^k_j)-h_k v_j^k
\end{eqnarray*}
and thus conclude that inclusion \eqref{evo_equa_discrete} holds at the discrete time $j$. It follows from the arguments below that for any $k\in\N$ sufficiently large we always have $\|x^k_{j+1}-\ox(t_{j+1}^k\|\le 1$. This completes the induction process and gives us therefore a sequence $\{(x^k(\cdot),u^k(\cdot))\}$ defined on the discrete mesh $T_k$ for large $k\in\N$ and satisfied therein the discretized sweeping inclusion \eqref{evo_equa_discrete} with the controlled moving set \eqref{mov-set}.

Next we prove that piecewise linear extensions $(x^k(t),u^k(t))$, $0\le t\le T$, of the above sequence to the continuous-time interval $[0,T]$ converges to the reference pair $(\ox(\cdot),\ou(t))$ in the norm topology of $W^{1,2}([0,T];\R^n)\times{\cal C}([0,T];\R^m)$. To proceed, fix any $\epsilon>0$ and recall the definition of $\epsilon_j$ in \eqref{ep}. Taking into account that $\ox(\cdot)\in{\cal C}^1([0,T];\R^n)$, we get that $\epsilon_j\le h_k\epsilon$ for all $j=0,\ldots,k-1$ and all $k\in\N$ sufficiently large. This gives us the relationships
\begin{align*}
\|x_{j+1}^k-\ox(t_{j+1}^k)\|&\le\|x^k_j+h_k f(t^k_j,x^k_j)-h_k v_j^k-\ox(t_{j+1}^k)\|\\
&\le\|x_j^k-\ox(t_j^k)\|+h_k\|f(t^k_j,x^k_j)-f(t^k_j,\ox(t^k_j))\|+h_k\|f(t^k_j,\ox(t^k_j))-v_j^k-\dot\ox(t_j^k)\|+\epsilon_j\\
&\le\big(1+(\alpha+L_f)h_k+\beta h_k\|x_j^k-\ox(t_j^k)\|\big)\|x_j^k-\ox(t_j^k)\|+\epsilon_j\\
&\le\big(1+(\alpha+L_f)h_k+\beta h_k\|x_j^k-\ox(t_j^k)\|\big)\|x_j^k-\ox(t_j^k)\|+h_k\epsilon,
\end{align*}
and therefore we arrive at the following estimate:
\begin{eqnarray*}
\|x_{j+1}^k-\ox(t_{j+1}^k)\|\le\big(1+(\alpha+L_f h_k+\beta h_k\|x_j^k-\ox(t_j^k)\|\big)\|x_j^k-\ox(t_j^k)\|+h_k\epsilon.
\end{eqnarray*}
Define further the quantities $a_j:=\|x_j^k-\ox(t_j^k)\|$ for all $j=0,\ldots,k$ and $k\in\N$. Observe that $a_{j+1}=(1+(\alpha+L_f) h_k+\beta h_k a_j)a_j+h_k \epsilon$ for $j=0,\ldots,k-1$, which is equivalent to
\begin{eqnarray*}
\frac{a_{j+1}-a_j}{h_k}=(\alpha+L_f)a_j+\beta a_j^2 +\epsilon,\quad j=0,\ldots,k-1.
\end{eqnarray*}
Denoting by $a_{\epsilon}(\cdot)$ the solution of the differential equation
\begin{eqnarray*}
\dot a(t)=(\alpha+L_f)a(t)+\beta a^2(t)+\epsilon,\quad t\in[0,T],\quad a(0)=0,
\end{eqnarray*}
we see that $a_\epsilon(t)\to 0$ uniformly on $[0,T]$ as $\epsilon\downarrow 0$. It readily implies that
\begin{eqnarray*}
\max\big\{\|x_j^k-\ox(t_j^k)\|,\;j=0,\ldots,k\big\}\to 0\;\mbox{ as }\;k\to\infty.
\end{eqnarray*}
This verifies, in particular, that $\max\{\|x_j^k-\ox(t_j^k)\|,\;j=0,\ldots,k\}\le 1$ for all $k\in\N$, which was needed to complete the induction process. Since we have
\begin{eqnarray}\label{est}
\|u_j^k-\ou(t_j^k)\|\le\tau L_g\|x_j^k-\ox(t_j^k)\|\;\mbox{ for all }\;j=0,\ldots,k
\end{eqnarray}
as shown above, the control sequence $\{u^k(t)\}$ converges to $\ou(\cdot)$ strongly in ${\cal C}([0,T];\R^m)$.

Let us next justify the strong $W^{1,2}$-convergence of the trajectories $x^k(\cdot)$ to $\ox(\cdot)$ on $[0,T]$. We have
\begin{align*}
\int_0^T \|\dot x^k(t)-\dot\ox(t)\|^2dt&=\sum_{j=0}^{k-1}\int_{t^k_j}^{t^k_{j+1}}\|f(t_j^k,x_j^k)-v_j^k-\dot\ox(t)\|^2dt\\
&\le 2\sum_{j=0}^{k-1}h_k\big\|f(t_j^k,x_j^k)-v_j^k-\dot\ox(t^k_j)\big\|^2+2\sum_{j=0}^{k-1}\int_{t^k_j}^{t^k_{j+1}}\|\dot\ox(t^k_j)-\dot\ox(t)\|^2dt,
\end{align*}
where the last term converges to zero due to $\ox(\cdot)\in{\cal C}^1([0,T];\R^n)$. The first term therein also converges to zero by the following estimates valid for all $j=0,\ldots,k$:
\begin{align*}
\|f(t_j^k,x_j^k)-v_j^k-\dot\ox(t^k_j)\|&\le\big\|f\big(t_j^k,\ox(t_j^k)\big)-v_j^k-\dot\ox(t^k_j)\big\|+\|f(t_j^k,x_j^k)-f(t_j^k,\ox(t_j^k))\|\\
&\le\big(L_f+\alpha+\beta\|x_j^k-\ox(t_j^k)\|\big)\|x_j^k-\ox(t^k_j)\|,
\end{align*}
which is due to \eqref{derivative_est}. Thus we get $\int_0^T\|\dot x^k(t)-\dot\ox(t)\|^2dt\to 0$ as $k\to\infty$. Since $x_0^k=\ox(0)$, the latter verifies that $\{x^k(\cdot)\}$ strongly converges to $\ox(\cdot)$ in $W^{1,2}([0,T];\R^n)$.

To complete the proof of the theorem, it remains to show that if $\ou(\cdot)\in W^{1,2}([0,T];\R^m)$ and the mapping $g(\cdot)$ is linear, then $u^k(\cdot)\to\ou(\cdot)$ strongly in $W^{1,2}([0,T];\R^m)$. Denoting
\begin{eqnarray*}
A^k_j:=\big(g(x_j^k),g(\ox(t_j^k)),\ou(t_j^k)\big)\;\mbox{ for all }\;j=0,\ldots,k-1\;\mbox{ and }\;k\in\N
\end{eqnarray*}
and using the local Lipschitz continuity of $d(\cdot)$ with constant $L_d>0$, we get
\begin{align*}
&\Big\|\frac{u^k_{j+1}-u_j^k}{h_k}-\frac{\ou(t_{j+1}^k)-\ou(t_j^k)}{h_k}\Big\|=\frac{1}{h_k}\big\|d\big(g(x_{j+1}^k)-g(\ox(t_{j+1}^k)),
A^k_{j+1}\big)-d\big(g(x_j^k)-g(\ox(t_j^k)),A^k_j\big)\big\|\\
&\le\frac{1}{h_k}\big\|d\big(g(x_{j+1}^k)-g(\ox(t_{j+1}^k)),A^k_{j+1}\big)-d\big(g(x_j^k)-g(\ox(t_j^k)),A^k_{j+1}\big)\big\|\\
&+\big\|d\big(g(x_j^k)-g(\ox(t_j^k)),A^k_{j+1}\big)-d\big(g(x_j^k)-g(\ox(t_j^k)),A^k_j\big)\big\|\\
&\le\frac{L_d L_g}{h_k}\big\|\big(x_{j+1}^k-\ox(t_{j+1}^k)\big)-\big(x_j^k-\ox(t_j^k)\big)\big\|\\
&+M\frac{L_d}{h_k}\big(\|g(x_{j+1}^k)-g(x_j^k)\|+\big\|g\big(\ox(t_{j+1}^k)\big)-g\big(\ox(t_j^k)\big)\big\|+\|\ou(t_{j+1}^k)-\ou(t_j^k)\|\big),
\end{align*}
where $M>0$ is sufficiently large. This shows that
\begin{align*}
&\int_0^T\|\dot u^k(t)-\dot\ou(t)\|^2dt\le M\int_0^T\|\dot x^k(t)-\dot\ox(t)\|^2dt\\
&+M\sum_{i=0}^{k-1}h_k\Big(\|g(x_{j+1}^k)-g(x_j^k)\|^2+\big\|g\big(\ox(t_{j+1}^k)\big)-g\big(\ox(t_j^k)\big)\big\|^2+\|\ou(t_{j+1}^k)-\ou(t_j^k)\|^2\Big)
\end{align*}
and thus verifies the claimed convergence under the assumptions made. $\h$\vspace*{0.05in}

The two approximation results established in Theorem~\ref{feasible-approx} allow us to apply the method of discrete approximations to deriving necessary optimality conditions for {\em two types} of local minimizers in problem $(P)$. The first type treats the trajectory and control components of the optimal pair $(\ox(\cdot),\ou(\cdot))$ in the same way and reduces in fact to the {\em intermediate $W^{1,2}$-minimizers} introduced in \cite{m95} in the general framework of differential inclusions and then studied in \cite{cm1}--\cite{h1} for various controlled sweeping processes. The second type seems to be {\em new in control theory}; it treats control and trajectory components differently and applies to problems $(P)$ whose running costs do not depend on control velocities.\vspace*{-0.05in}

\begin{Definition}{\bf(local minimizers for controlled sweeping processes).}\label{ilm} Let the pair $(\ox(\cdot),\ou(\cdot))$ be feasible to problem $(P)$ under the standing assumptions made.

{\bf(i)} We say that $(\ox(\cdot),\ou(\cdot))$ be a {\sc local $W^{1,2}\times W^{1,2}$-minimizer} for $(P)$ if $\ox(\cdot)\in W^{1,2}([0,T];\R^n)$, $\ou(\cdot)\in W^{1,2}([0,T];\R^m)$, and
\begin{eqnarray}\label{lm1}
J[\ox,\ou]\le J[x,u]\;\mbox{ for all }\;x(\cdot)\in W^{1,2}([0,T];\R^n)\;\mbox{ and }\;u(\cdot)\in W^{1,2}([0,T];\R^m)
\end{eqnarray}
sufficiently close to $(\ox(\cdot),\ou(\cdot))$ in the norm topology of the corresponding spaces in \eqref{lm1}.

{\bf(ii)} Let the running cost $\ell(\cdot)$ in \eqref{eq:MP} do not depend on $\dot u$. We say that the pair $(\ox(\cdot),\ou(\cdot))$ be a {\sc local $W^{1,2}\times{\cal C}$-minimizer} for $(P)$ if $\ox(\cdot)\in W^{1,2}([0,T];\R^n)$, $\ou(\cdot)\in{\cal C}([0,T];\R^m)$, and
\begin{eqnarray}\label{lm2}
J[\ox,\ou]\le J[x,u]\;\mbox{ for all }\;x(\cdot)\in W^{1,2}([0,T];\R^n)\;\mbox{ and }\;u(\cdot)\in{\cal C}([0,T];\R^m)
\end{eqnarray}
sufficiently close to $(\ox(\cdot),\ou(\cdot))$ in the norm topology of the corresponding spaces in \eqref{lm2}.
\end{Definition}\vspace*{-0.05in}

Our main attention in what follows is to derive {\em necessary optimality conditions} for both types of local minimizers in Definition~\ref{ilm} by developing appropriate versions of the method of discrete approximations. It is clear that any local $W^{1,2}\times{\cal C}$-minimizer for $(P)$ is also a local $W^{1,2}\times W^{1,2}$-minimizer for this problem, provided that we restrict the class of feasible controls to $W^{1,2}$-functions. Thus necessary optimality conditions for local $W^{1,2}\times W^{1,2}$-minimizers are also necessary for local $W^{1,2}\times{\cal C}$-ones in this framework, while not vice versa. On the other hand, we may deal with local $W^{1,2}\times{\cal C}$-minimizers without imposing anything but the continuity assumptions of feasible controls, provided that the running cost in \eqref{eq:MP} does not depend on control velocities. Note furthermore that considering a $W^{1,2}$-neighborhood of the trajectory part $\ox(\cdot)$ in both settings of Definition~\ref{ilm} leads us to potentially more selective necessary optimality conditions for such minimizers than for conventional strong local minimizers and global solutions to $(P)$.

It has been well recognized in the calculus of variations and optimal control, starting with pioneering studies by Bogolyubov and Young, that limiting procedures of dealing with continuous-time dynamical systems involving time derivatives require a certain {\em relaxation stability}, which means that the value of cost functionals does not change under the convexification of the dynamics and running cost with respect to velocity variables; see, e.g., \cite{mord,vinter} for more details and references. In sweeping control theory, such issues have been investigated in \cite{et,t} for controlled sweeping processes somewhat different from $(P)$.

To consider an appropriate relaxation of our problem $(P)$, denote
\begin{eqnarray}\label{F}
F=F(t,x,u):=f(t,x)-N\big(g(x);C(t,u)\big)
\end{eqnarray}
and formulate the {\em relaxed optimal control problem} $(R)$ as a counterpart of $(P)$ with the replacement of the cost functional \eqref{eq:MP} by the convexified one
\begin{eqnarray*}\label{R}
{\rm minimize}\;\Hat J[x,u]:=\varphi\big(x(T)\big)+\int_0^T\Hat\ell_F\big(t,x(t),u(t),\dot x(t),\dot u(t)\big)dt,
\end{eqnarray*}
where $\Hat\ell(t,x,u,\cdot,\cdot)$ is defined as the largest l.s.c.\ convex function majorized by $\ell(t,x,u,\cdot,\cdot)$ on the convex closure of the set $F$ in \eqref{F} with $\Hat\ell:=\infty$ otherwise. Then we say that the pair $(\ox(\cdot),\ou(\cdot))$ is a {\em relaxed local $W^{1,2}\times W^{1,2}$-minimizer} for $(P)$ if in additions to the conditions of Definition~\ref{ilm}(i) we have $J[\ox,\ou)=\Hat J[\ox,\ou]$. Similarly we define a {\em relaxed local $W^{1,2}\times{\cal C}$-minimizer} for $(P)$ in the setting of Definition~\ref{ilm}(ii). Note that, in contrast to the original problem $(P)$, the convexified structure of the relaxed problem $(R)$ provides an opportunity to the establish the {\em existence} of global optimal solutions in the prescribed classes of controls and trajectories. It is not a goal of this paper, but we refer the reader to \cite[Theorem~4.1]{cm3} and \cite[Theorem~4.2]{t} for some particular settings of controlled sweeping processes in the classes of $W^{1,2}\times W^{1,2}$ and $W^{1,2}\times{\cal C}$ feasible pairs $(\ox(\cdot),\ou(\cdot))$, respectively.

There is clearly no difference between the problems $(P)$ and $(R)$ if the normal cone in \eqref{F} is convex and the integrand $\ell$ in \eqref{eq:MP} is convex with respect to velocity variables. On the other hand, the measure continuity/nonatonomicity on $[0,T]$ and the differential inclusion structure of the sweeping process \eqref{evo_equa} create the environment where any local minimizer of the types under consideration is also a relaxed one. Without delving into details here, we just mention that the possibility to derive such a {\em local relaxation stability} from \cite[Theorem~4.2]{t} for {\em strong} local (in the ${\cal C}$-norm) minimizers of $(P)$, provided that the controlled moving set $C(t,u)$ in \eqref{mov-set} is convex and continuously depends on its variables.\vspace*{0.02in}

Given now a relaxed local minimizer $(\ox(\cdot),\ou(\cdot))$ of the types introduced in Definition~\ref{ilm}, we construct appropriate sequences of discrete-time optimal control problems corresponding to each type therein separately. For {\em brevity and simplicity}, from now on we restrict ourselves to the setting of $(P)$ where $g(x):=x$, $f:=0$ while $\psi$ and $\ell$ do not depend on $t$. The reader can easily check that the procedure developed below is applicable to the general version of $(P)$ under the standing assumptions made.

If the pair $(\ox(\cdot),\ou(\cdot))$ is a relaxed local {\em $W^{1,2}\times W^{1,2}$-minimizer} of $(P)$, we fix $\ve>0$ sufficiently small to accommodate the $W^{1,2}\times W^{1,2}$-neighborhood of $(\ox(\cdot),\ou(\cdot))$ in Definition~\ref{ilm}(i) and for each $k\in\N$ define the approximation problem $(P^1_k)$ as follows:
\begin{eqnarray*}\label{cost-pk}
\begin{array}{ll}
\mbox{minimize}\;\;
J_k[z^k]:=&\disp\varphi(x_k^k)+\disp{h_k}\sum\limits_{j=0}^{k-1}\disp{\ell\Big(x_j^k,u_j^k,\frac{x_{j+1}^k-x_j^k}{h_k},
\disp\frac{u_{j+1}^k-u_j^k}{h_k}\Big)}\\
&+h_k\disp\sum\limits_{j=0}^{k-1}\int\limits_{{t^k_j}}^{{t^k_{j+1}}}\disp\Big(\Big\|\frac{x_{j+1}^k-x_j^k}{h_k}-\dot\ox(t)\Big\|^2+
\Big\|\disp\frac{u_{j+1}^k-u_j^k}{h_k}-\dot{\ou}(t)\Big\|^2\Big)\,dt
\end{array}
\end{eqnarray*}
over collections $z^k:=(x_0^k,\ldots,x_k^k,u_0^k,\ldots,u_k^k)$ subject to the constraints
\begin{eqnarray}\label{discrete-inclusion}
x_{j+1}^k\in x_j^k+h_k F(x_j^k,u_j^k)\;\mbox{ for }\;j=0,\ldots,k-1\;\mbox{ with }\;\big(x^k_0,u^k_0\big)=\big(x_0,\ou(0)\big),
\end{eqnarray}
\begin{eqnarray}\label{state-constraint}
(x_k^k,u^k_k)\in\psi^{-1}(\Th),
\end{eqnarray}
\begin{eqnarray}\label{nei}
\big\|(x^k_j,u^k_j)-\big(\ox(t^k_j),\ou(t^k_j)\big)\big\|\le\ve/2\;\mbox{ for }\;j=0,\ldots,k,
\end{eqnarray}
\begin{eqnarray*}\label{neighborhood1}
\sum\limits_{j=0}^{k-1}\int\limits_{{t^k_j}}^{{t^k_{j+1}}}\disp\Big(\Big\|\frac{x_{j+1}^k-x_j^k}{h_k}\disp-
\dot\ox(t)\Big\|^2+\Big\|\disp\frac{u_{j+1}^k-u_j^k}{h_k}-\dot\ou(t)\Big\|^2\Big)\,dt\le\frac{\ve}{2}.
\end{eqnarray*}

If the pair $(\ox(\cdot),\ou(\cdot))$ is a relaxed local {\em $W^{1,2}\times{\cal C}$-minimizer} of $(P)$, fix $\ve>0$ sufficiently small to accommodate the $W^{1,2}\times{\cal C}$-neighborhood of $(\ox(\cdot),\ou(\cdot))$ in Definition~\ref{ilm}(ii) and for each $k\in\N$ define the approximation problem $(P^2_k)$ in the following way corresponding to \eqref{lm2}:
\begin{eqnarray*}\label{cost-pk1}
\begin{array}{ll}
\mbox{minimize}\;\;
J_k[z^k]:=&\disp\varphi(x_k^k)+\disp{h_k}\sum\limits_{j=0}^{k-1}\disp{\ell\Big(x_j^k,u_j^k,\frac{x_{j+1}^k-x_j^k}{h_k}\Big)}\\
&+\disp\sum\limits_{j=0}^{k}\big\|u_j^k-\ou(t_j^k)\big\|^2+\disp\sum\limits_{j=0}^{k-1}\int\limits_{{t^k_j}}^{{t^k_{j+1}}}\disp\Big\|\frac{x_{j+1}^k-x_j^k}{h_k}
-\dot\ox(t)\Big\|^2dt
\end{array}
\end{eqnarray*}
over $z^k=(x_0^k,\ldots,x_k^k,u_0^k,\ldots,u_k^k)$  subject to the constraints in \eqref{discrete-inclusion}--\eqref{nei} and
\begin{eqnarray*}\label{neighborhood2}
\sum\limits_{j=0}^{k-1}\int\limits_{{t^k_j}}^{{t^k_{j+1}}}\disp\Big\|\frac{x_{j+1}^k-x_j^k}{h_k}\disp-
\dot\ox(t)\Big\|^2 dt\le\frac{\ve}{2}.
\end{eqnarray*}

To proceed further with the method of discrete approximations, we need to make sure that the approximating problems $(P^i_k)$, $i=1,2$, admit optimal solutions. This is indeed the case due to Theorem~\ref{feasible-approx} and the {\em robustness} (closed-graph property) of our basic normal cone \eqref{nc}.\vspace*{-0.05in}

\begin{Proposition}{\bf (existence of discrete optimal solutions).}\label{ex-disc} Under the imposed standing assumptions {\rm(H1)--(H5)}, each problem $(P^i_k)$, $i=1,2$, has an optimal solution for all $k\in\N$ sufficiently large.
\end{Proposition}\vspace*{-0.05in}
{\bf Proof.} It follows from Theorem~\ref{feasible-approx} and the constructions of $(P^1_k)$, $(P^2_k)$ that the set of feasible solutions to each of these problems is nonempty for all large $k\in\N$. Applying the classical Weierstrass existence theorem, observe that the boundedness of the feasible sets follows directly from the constraint structures in $(P^1_k)$ and $(P^2_k)$. The remaining closedness of the feasible sets for these problems is a consequence of the robustness property of the normal cone \eqref{nc} that determines the discrete inclusions \eqref{discrete-inclusion}. $\h$\vspace*{0.05in}

The next theorem establishes the strong convergence in the corresponding spaces of extended discrete optimal solutions for discrete approximation problems to the given relaxed local minimizers for $(P)$.\vspace*{-0.05in}

\begin{Theorem}{\bf(strong convergence of discrete optimal solutions).}\label{conver} In addition to the standing assumption {\rm(H1)--(H5)}, suppose that the cost functions $\ph$ and $\ell$ are continuous around the given local minimizer. The following assertions hold:

{\bf(i)} If $(\ox(\cdot),\ou(\cdot))$ is a relaxed local $W^{1,2}\times W^{1,2}$-minimizer for $(P)$, then any sequence of piecewise linear extensions on $[0,T]$ of the optimal solutions $(\ox^k(\cdot),\ou^k(\cdot))$ to $(P^1_k)$ converges to $(\ox(\cdot),\ou(\cdot))$ in the norm topology of $W^{1,2}([0,T];\R^n)\times  W^{1,2}([0,T];\R^m)$ as $k\to\infty$.

{\bf(ii)} If $(\ox(\cdot),\ou(\cdot))$ is a relaxed local $W^{1,2}\times{\cal C}$-minimizer for $(P)$, then any sequence of piecewise linear extensions on $[0,T]$ of the optimal solutions $(\ox^k(\cdot),\ou^k(\cdot))$ to $(P^2_k)$ converges to $(\ox(\cdot),\ou(\cdot))$ in the norm topology of $W^{1,2}([0,T];\R^n)\times{\cal C}([0,T];\R^m)$ as $k\to\infty$.
\end{Theorem}\vspace*{-0.05in}
{\bf Proof.} To verify assertion (i), we proceed similarly to the proof of \cite[Theorem~3.4]{h1} with the usage of the normal cone robustness for \eqref{nc} instead of Attouch's theorem in the convex setting of \cite{h1}. The proof of assertion (ii) goes in the same lines with observing that the required ${\cal C}$-compactness of the control sequence follows from the control-state relationship of type \eqref{est} valid due to assumption (H4). $\h$\vspace*{-0.1in}

\section{Extended Euler-Lagrange Conditions for Sweeping Solutions}\label{optimality-conditions}
\setcounter{equation}{0}

Having the strong convergence results of Theorem~\ref{conver} as the quintessence of the discrete approximation well-posedness justified in Section~\ref{discrete-approximations}, we now proceed first with deriving necessary optimality conditions in both discrete problems $(P^1_k)$ and $(P^2_k)$ for each $k\in\N$ and then with the subsequent passage to the limit therein as $k\to\infty$. In this way we arrive at necessary optimality conditions for relaxed local minimizers in $(P)$ of both $W^{1,2}\times W^{1,2}$ and $W^{1,2}\times{\cal C}$ types.

Observe that for each fixed $k\in\N$ both problems $(P^1_k)$ and $(P^2_k)$ belong to the class of finite-dimensional mathematical programming with nonstandard geometric constraints \eqref{discrete-inclusion} and \eqref{state-constraint}. We can handle them by employing appropriate tools of variational analysis that revolve around the normal cone \eqref{nc}.\vspace*{-0.05in}

\begin{Theorem}{\bf(necessary optimality conditions for $(P^1_k)$).}\label{discrete} Fix $k\in\N$ and consider an optimal solution $\oz^k:=(x_0,\ox^k_1\ldots,\ox_k^k,\ou_0^k,\ldots,\ou_{k}^{k})$ to problem $(P^1_k)$, where $F$ may be a general closed-graph mapping. Suppose that the cost functions $\ph$ and $\ell$ are locally Lipschitzian around the corresponding components of the optimal solution and denote the quantities
\begin{eqnarray}\label{th}
\Big(\theta_{j}^{xk},\theta_{j}^{uk}\Big):=2\int\limits_{{t_j^k}}^{{t_{j+1}^k}}\Big(\frac{\ox_{j+1}^k-\ox_j^k}{h_k}
\disp-\dot\ox(t),\frac{\ou_{j+1}^k-\ou_j^k}{h_k}-\dot\ou(t)\Big)dt,\quad j=0,\ldots,k-1.
\end{eqnarray}
Then there exist dual elements $\lm^k\ge 0$, $p_j^k=(p^{xk}_j,p^{uk}_j)\in\R^n\times\R^m$ as $j=0,\ldots,k$ and subgradient vectors
\begin{eqnarray}\label{sublneighborhood1}
\big(w^{xk}_{j},w^{uk}_{j},v^{xk}_{j},v^{uk}_{j}\big)\in\partial\ell\left(\ox^k_j,\ox^k_j,\frac{\ox_{j+1}^k-\ox_j^k}{h_k},\frac{\ou_{j+1}^k-\ou_j^k}{h_k}\right),\quad j=0,\ldots,k-1,
\end{eqnarray}
such that the following conditions are satisfied:
\begin{eqnarray}\label{nontriv}
\lm^k+\sum_{j=0}^{k-1}\|p^{xk}_j\|+\|p^{uk}_0\|+\|p^{xk}_k\|+\|p^{uk}_k\|\ne 0,
\end{eqnarray}
\begin{eqnarray}\label{transversality_end_d}
-(p^{xk}_k,p^{uk}_k)\in\lambda^k\big(\partial\varphi(\ox_k^k),0\big)+N\big((\ox^k_k,\ou^k_k);\psi^{-1}(\Th)\big),
\end{eqnarray}
\begin{eqnarray}\label{psiu}
p^{uk}_{j+1}=\lambda^k(v^{uk}_j+h_k^{-1}\theta^{uk}_j),\;j=0,\ldots,k-1,
\end{eqnarray}
\begin{eqnarray}\label{euler1}
\begin{split}
&\left(\frac{p_{j+1}^{xk}-p_j^{xk}}{h_k}-\lambda^k w_j^{xk},\frac{p_{j+1}^{uk}-p_j^{uk}}{h_k}-
\lambda^k w_j^{uk},p_{j+1}^{xk}-\lambda^k\Big(v_j^{xk}+\frac{1}{h_k}\theta_j^{xk}\Big)\right)\\
&\qquad\in N\Big(\Big(\ox_j^k,\ou_j^k,\frac{\ox_{j+1}^k-\ox_j^k}{h_k}\Big);\gph F\Big),\quad j=0,\ldots,k-1.
\end{split}
\end{eqnarray}
\end{Theorem}\vspace*{-0.05in}
{\bf Proof.} It follows the lines in the proof of \cite[Theorem~5.1]{h1} by reducing $(P^1_k)$ to a problem of mathematical programming. The usage of necessary optimality conditions for such problems and calculus rules of generalized differentiation for the basic constructions \eqref{nc} and \eqref{1sub} available in the books \cite{mordukhovich,mord,rw} allow us arrive at \eqref{sublneighborhood1}--\eqref{euler1} due to the particular structure of the data in $(P^1_k)$. $\h$\vspace*{0.05in}

The same approach holds for verifying the necessary optimality conditions for problem $(P^2_k)$ presented in the next theorem, which also takes into account the specific structure of this problem.\vspace*{-0.07in}

\begin{Theorem}{\bf(necessary optimality conditions for $(P^2_k)$).}\label{discrete1} Let $\oz^k:=(x_0,\ox^k_1\ldots,\ox_k^k,\ou_0^k,\ldots,\ou_{k}^{k})$ be an optimal solution problem $(P^2_k)$ in the framework of Theorem~{\rm\ref{discrete}}. Consider the quantities
\begin{eqnarray*}
\theta_{j}^{xk}:=2\int\limits_{{t_j^k}}^{{t_{j+1}^k}}\Big(\frac{\ox_{j+1}^k-\ox_j^k}{h_k}
\disp-\dot\ox(t)\Big)dt,\quad\theta_{j}^{uk}:=2\big(\ou_j^k-\ou(t_j^k)\big),\quad j=0,\ldots,k.
\end{eqnarray*}
Then there exist dual elements $\lm^k\ge 0$, $p_j^k=(p^{xk}_j,p^{uk}_j)\in\R^n\times\R^m$ as $j=0,\ldots,k$ and subgradient vectors
\begin{eqnarray*}
\big(w^{xk}_{j},w^{uk}_{j},v^{xk}_{j}\big)\in\partial\ell\left(\ox^k_j,\ou^k_j,\frac{\ox_{j+1}^k-\ox_j^k}{h_k}\right),\quad j=0,\ldots,k-1,
\end{eqnarray*}
satisfying following necessary optimality conditions:
\begin{eqnarray*}
\lm^k+\sum_{j=0}^{k-1}\|p^{xk}_j\|+\|p^{uk}_0\|+\|p^{xk}_k\|\ne 0,
\end{eqnarray*}
\begin{eqnarray*}
-(p^{xk}_k,0)\in\lambda^k\big(\partial\varphi(\ox_k^k),0\big)+N\big((\ox^k_k,\ou^k_k);\psi^{-1}(\Th)\big),
\end{eqnarray*}
\begin{eqnarray*}
p^{uk}_{j+1}=0,\quad j=0,\ldots,k-1,
\end{eqnarray*}
\begin{eqnarray*}
\begin{split}
&\left(\frac{p_{j+1}^{xk}-p_j^{xk}}{h_k}-\lambda^k w_j^{xk},\frac{p_{j+1}^{uk}-p_j^{uk}}{h_k}-
\lambda^k (w_j^{uk}+\theta_{j}^{uk}),p_{j+1}^{xk}-\lambda^k\Big(v_j^{xk}+\frac{1}{h_k}\theta_j^{xk}\Big)\right)\\
&\qquad\in N\Big(\Big(\ox_j^k,\ou_j^k,\frac{\ox_{j+1}^k-\ox_j^k}{h_k}\Big);\gph F\Big),\quad j=0,\ldots,k-1.
\end{split}
\end{eqnarray*}
\end{Theorem}\vspace*{-0.05in}

Now we are ready to derive necessary optimality conditions for both types of (relaxed) local minimizers for $(P)$ from Definition~\ref{ilm} by passing the limit from those in Theorems~\ref{discrete} and \ref{discrete1} with taking into account the convergence results from Section~\ref{discrete-approximations} and calculation results of generalized differentiation presented in Section~\ref{prelim}. The reader can see that the obtained optimality conditions for both types are local minimizers are pretty similar under the imposed assumptions. This is largely due to the achieved discrete approximation convergence in Theorem~\ref{conver} and the structures of the discretized problems. Necessary optimality conditions for relaxed local $W^{1,2}\times W^{1,2}$-minimizers of $(P)$ were derived in \cite[Theorem~6.1]{h1} for polyhedral moving sets \eqref{sw-con2} under significantly more restrictive assumptions, which basically cover the case of $(\ox(\cdot),\ou(\cdot))\in W^{2,\infty}([0,T];\R^{n+m})$. Note that the {\em linear independence constraint qualification} (LICQ) condition on generating polyhedral vectors imposed therein is a counterpart of our surjectivity assumption (H3) in the polyhedral setting of \cite{h1}.\vspace*{-0.05in}

\begin{Theorem}{\bf(necessary optimality conditions for the controlled sweeping process).}\label{necopt} Let $(\ox(\cdot),\ou(\cdot))$ be a local minimizer for problem $(P)$ of the types specified below. In addition to the standing assumptions, suppose that $\psi=\psi(x,u)$ is ${\cal C}^2$-smooth with respect to both variables while $\ph$ and $\ell$ are locally Lipschitzian around the corresponding components of the optimal solution. The following assertions hold:

{\bf(i)} If $(\ox(\cdot),\ou(\cdot))$ is a relaxed local $W^{1,2}\times W^{1,2}$-minimizer, then there exist a multiplier $\lm\ge 0$, an adjoint arc $p(\cdot)=(p^x,p^u)\in W^{1,2}([0,T];\R^n\times\R^m)$, a signed vector measure $\gamma\in C^*([0,T];\R^s)$, as well as pairs
$(w^x(\cdot),w^u(\cdot))\in L^2([0,T];\R^n\times\R^m)$ and $(v^x(\cdot),v^u(\cdot))\in L^\infty([0,T];\R^n\times\R^m)$ with
\begin{eqnarray}\label{co}
\big(w^x(t),w^u(t),v^x(t),v^u(t)\big)\in{\rm co}\,\partial\ell\big(\ox(t),\ou(t),\dot\ox(t),\dot\ou(t)\big)\;\mbox{ a.e. }\;t\in[0,T]
\end{eqnarray}
satisfying the collection of necessary optimality conditions:

$\bullet$ {\sc Primal-dual dynamic relationships}:
\begin{eqnarray}\label{hamiltonx}
\dot p(t)=\lm w(t)+\left[\begin{array}{c}
\nabla_{xx}^{2}\big\langle\eta(t),\psi\big\rangle\big(\ox(t),\ou(t)\big)\\
\nabla_{xw}^{2}\big\langle\eta(t),\psi\big\rangle\big(\ox(t),\ou(t)\big)
\end{array}\right]\big(-\lm v^x(t)+q^x(t)\big)\;\mbox{ a.e. }\;t\in[0,T],
\end{eqnarray}
\begin{eqnarray}\label{hamilton2ub}
q^u(t)=\lm v^u(t)\;\mbox{ a.e. }\;t\in[0,T],
\end{eqnarray}
where $\eta(\cdot)\in L^{2}([0,T];\R^s)$ is a uniquely defined vector function determined by the representation
\begin{eqnarray}\label{etajkl}
\dot\ox(t)=-\nabla_x\psi\big(\ox(t),\ou(t)\big)^*\eta(t)\;\mbox{ a.e. }\;t\in[0,T]
\end{eqnarray}
with $\eta(t)\in N(\psi(\ox(t),\ou(t));\Th)$, and where $q\colon[0,T]\to\R^n\times\R^m$ is a function of bounded variation on $[0,T]$ with its left-continuous
representative given, for all $t\in[0,T]$ except at most a countable subset, by
\begin{eqnarray}\label{q}
q(t)=p(t)-\int_{[t,T]}\nabla\psi\big(\ox(\t),\ou(\t)\big)^*d\gamma(\t).
\end{eqnarray}

$\bullet$ {\sc Measured coderivative condition}: Considering the $t$-dependent outer limit
\begin{eqnarray*}
\Limsup_{|B|\to 0}\frac{\gamma(B)}{|B|}(t):=\Big\{y\in\R^s\Big|\;\exists\,\mbox{ sequence }\;B_k\subset[0,1]\;\mbox{ with }\;t\in\B_k,\;|B_k|\to 0,\;\frac{\gamma(B_k)}{|B_k|}\to y\Big\}
\end{eqnarray*}
over Borel subsets $B\subset[0,1]$ with the Lebesgue measure $|B|$, for a.e.\ $t\in[0,T]$ we have
\begin{eqnarray}\label{maximumcondition}
D^*N_\Th\big(\psi(\ox(t),\ou(t)),\eta(t)\big)\big(\nabla_x\psi(\ox(t),\ou(t))(q^x(t)-\lm v^x(t))\big)\cap\Limsup_{|B|\to 0}\frac{\gamma(B)}{|B|}(t)\ne\emp.
\end{eqnarray}

$\bullet$ {\sc Transversality condition} at the right endpoint:
\begin{eqnarray}\label{pxkkc}
-\big(p^x(T),p^u(T)\big)\in\lm\big(\partial\ph(\ox(T)),0\big)+\nabla\psi\big(\ox(T),\ou(T)\big)N_\Th\big((\ox(T),\ou(T)\big).
\end{eqnarray}

$\bullet$ {\sc Measure nonatomicity condition}: Whenever $t\in[0,T)$ with $\psi(\ox(t),\ou(t))\in{\rm int}\,\Th$ there is a neighborhood $V_t$ of $t$ in $[0,T]$ such that $\gamma(V)=0$ for any Borel subset $V$ of $V_t$.

$\bullet$ {\sc Nontriviality condition}:
\begin{eqnarray}\label{nontriv2}
\lm+\sup_{t \in[0,T]}\|p(t)\|+\|\gamma\|\ne 0\;\mbox{ with }\;\|\gamma\|:=\sup_{\|x\|_{C([0,T]}=1}\int_{[0,T]}x(s)d\gamma.
\end{eqnarray}

{\bf(ii)} If $(\ox(\cdot),\ou(\cdot))$ is a relaxed local $W^{1,2}\times{\cal C}$-minimizer, then all the conditions \eqref{hamiltonx}--\eqref{nontriv2} in {\rm(i)} hold with the replacement of the quadruple $(w^x(\cdot),w^u(\cdot),v^x(\cdot),v^u(\cdot))$ in \eqref{co} by the triple $(w^x(\cdot),w^u(\cdot),v^x(\cdot))\in L^2([0,T];\R^n)\times L^2([0,T];\R^m)\times L^\infty([0,T];\R^n)$ satisfying the inclusion
\begin{eqnarray*}
\big(w^x(t),w^u(t),v^x(t)\big)\in{\rm co}\,\partial\ell\big(\ox(t),\ou(t),\dot\ox(t)\big)\;\mbox{ a.e. }\;t\in[0,T].
\end{eqnarray*}
\end{Theorem}\vspace*{-0.05in}
{\bf Proof.} We give it only for assertion (i), since the proof of (ii) is similar with taking into account the type of convergence $\ou^k(\cdot)\to\ou(\cdot)$ achieved in Theorem~\ref{conver}(ii) and that the running cost $\ell$ in Definition~\ref{ilm}(ii) does not depend on the control velocity $\dot{u}$.

To verify assertion (i), deduce first from \eqref{euler1} and Proposition~\ref{morout} that for each $k\in\N$ and $j=0,\ldots,k-1$ there is a unique vector $\eta_j^k\in N_{Z}(\psi(\ox_j^k,\ou_j^k))$ satisfying the conditions
\begin{eqnarray*}\label{eta_dis}
\nabla_{x}\psi(\ox_j^k,\ou_j^k)^*\eta_j^k=-\frac{\ox_{j+1}^k-\ox_j^k}{h_k},
\end{eqnarray*}
\begin{eqnarray}\label{p_iteration}
\frac{p_{j+1}^{xk}-p_j^{xk}}{h_k}-\lambda^k w_j^{k}=\left[\begin{array}{c}
\nabla_{xx}^{2}\langle\eta_j^k,\psi\rangle(\ox_j^k,\ou_j^k)\\
\nabla_{xw}^{2}\langle\eta_j^k,\psi\rangle(\ox_j^k,\ou_j^k)
\end{array}\right]u+\nabla\psi(\ox_j^k,\ou_j^k)^*\gamma_j^k
\end{eqnarray}
with $u:=p_{j+1}^{xk}-\lambda^k\Big(v_j^{xk}+\frac{1}{h_k}\theta_j^{xk}\Big)$ and some vectors
\begin{eqnarray}\label{gamma_dis}
\gamma_j^k\in D^*N_{Z}\big(\psi(\ox_j^k,\ou_j^k),\eta_j^k\big)\Big(\nabla_{x}\psi(\ox_j^k,\ou_j^k)\Big(p_{j+1}^{xk}-\lambda^k\Big(v_j^{xk}+\frac{1}{h_k}\theta_j^{xk}\Big)\Big)
\Big).
\end{eqnarray}
Taking this into account, we get from \eqref{nontriv} the improved nontriviality condition
\begin{eqnarray}\label{nontriv1}
\lm^k+\|p^{uk}_0\|+\|p^{xk}_k\|+\|p^{uk}_k\|\ne 0\;\mbox{ for all }\;k\in\N
\end{eqnarray}
with the validity of \eqref{psiu} as well as $\lm^k\ge 0$ and the relationships in \eqref{th} and \eqref{sublneighborhood1} of Theorem~\ref{discrete}.

Now we proceed with passing to the limit as $k\to\infty$ in the obtained optimality conditions for discrete approximations. Since some arguments in this procedure are similar to those used in \cite[Theorem~6.1]{h1} in a more special setting, we skip them for brevity while focusing on significantly new developments. In particular, the existence of the claimed quadruples $(w^x(\cdot),w^u(\cdot),v^x(\cdot),v^u(\cdot))$ satisfying \eqref{co} is proved as in \cite{h2} while the existence of the uniquely defined $\eta(\cdot)\in L^2([0.T];\Th)$ solving the differential equation \eqref{etajkl} follows from representation \eqref{unique_re} by repeating the limiting procedure of \cite[Theorem~6.1]{h1}.

Next we define $q^k(\cdot)=(q^{xk}(\cdot),q^{uk}(\cdot))$ by extending $p_j^k$ piecewise linearly to $[0,T]$ with $q^{k}(t^k_j):=p^{k}_j$ for $j=0,\ldots,k$.
Construct $\gg^k(\cdot)$ on $[0,T]$ by
\begin{eqnarray*}
\gamma^k(t):=\gamma^k_j,\;\mbox{ for }\;t\in[t^k_j,t^k_{j+1}),\quad j=0,\ldots,k-1,
\end{eqnarray*}
with $\gamma^k(T):=0$ and consider the auxiliary functions
\begin{eqnarray}\label{vt}
\vt^{k}(t):=\max\big\{t_{j}^{k}\big|\;t_{j}^{k}\le t,\;0\le j\le k\big\}\;\mbox{ for all }\;
t\in\left[0,T\right],\;k\in\N,
\end{eqnarray}
so that $\vt^k(t)\to t$ uniformly in $[0,T]$ as $k\to\infty$. Since $\vt^{k}(t)=t_{j}^{k}$ for all $t\in\lbrack t_{j}^{k},t_{j+1}^{k})$ and $j=0,\ldots,k-1$, the equations in \eqref{p_iteration} can be rewritten as
\begin{eqnarray}\label{qk}
\dot q^{k}(t)-\lambda^k w^{k}(t)=\left[\begin{array}{c}
\nabla_{xx}^{2}\big\langle\eta^k_i(t),\psi\big\rangle\big(\ox^k(\vt^k(t),\ou^k(\vt^k(t)\big)\\
\nabla_{xw}^{2}\big\langle\eta^k_i(t),\psi\big\rangle\big(\ox^k(\vt^k(t),\ou^k(\vt^k(t)\big)
\end{array}\right]u+\nabla\psi\big(\ox^k(\vt^k(t),\ou^k(\vt^k(t)\big)^*\gamma^k(t),
\end{eqnarray}
where $u:=q^{xk}(\vt^k_+(t))-\lambda^k(v^{xk}(t)+\th^{xk}(t))$ for every $t\in(t^k_j,t^k_{j+1})$, $j=0,\ldots,k-1$, and $i=1,\ldots,m$, and where $\vt^k_+(t):=t^k_{j+1}$ for $t\in[t^k_j,t^k_{j+1})$.

Define now we $p^k(\cdot)=(p^{xk}(\cdot),p^{uk}(\cdot))$ on $[0,T]$ by setting
\begin{eqnarray*}
\begin{array}{ll}
p^{k}(t):=q^{k}(t)+\disp\int_{t}^{T}\nabla\psi\big(\ox^k(\vt^k(\t)),\ou^k(\vt^k(\t))\big)\gamma^k(\t)d\t
\end{array}
\end{eqnarray*}
for every $t\in[0,T]$. This gives us $p^k(T)=q^k(T)$ with the differential relation
\begin{eqnarray*}
\dot p^{k}(t)=\dot q^{k}(t)-\nabla\psi\big(\ox^k(\vt^k(t)),\ou^k(\vt^k(t))\big)^*\gamma^k(t)
\end{eqnarray*}
holding for a.e.\ $t\in[0,T]$. Substituting the latter into \eqref{qk}, we get
\begin{eqnarray*}
\dot p^{k}(t)-\lambda^k w^{k}(t)=\left[\begin{array}{c}
\nabla_{xx}^{2}\big\langle\eta^k_i(t),\psi\big\rangle\big(\ox^k(\vt^k(t)),\ou^k(\vt^k(t)\big)\\
\nabla_{xw}^{2}\big\langle\eta^k_i(t),\psi\big\rangle\big(\ox^k(\vt^k(t)),\ou^k(\vt^k(t)\big)
\end{array}\right]u,
\end{eqnarray*}
for every $t\in(t^k_j,t^k_{j+1})$, $j=0,\ldots,k-1$, and $i=1,\ldots,m$. Define further the vector measures $\gamma^k$ by
\begin{eqnarray*}
\int_B d\gamma^k:=\int_B\gamma^k(t)dt \;\mbox{ for every Borel subset }\;B\subset[0,T]
\end{eqnarray*}
and observe that, due to the  positive homogeneity of all the expressions in the statement of Theorem~\ref{discrete} with respect to $(\lambda^k,p^k,\gamma^k)$, the nontriviality condition \eqref{nontriv1} can be rewritten as
\begin{eqnarray}\label{alphak}
\lm^k+\|q^{uk}(0)\|+\|p^k(T)\|+\int_0^T\|\gamma^{k}(t)\|dt=1\;\mbox{ for all }\;k\in\N,
\end{eqnarray}
which tell us that all the sequential terms in \eqref{alphak} are uniformly bounded. Following the proof of \cite[Theorem 6.1]{h1}, we obtain the relationships in \eqref{hamiltonx}, \eqref{q}, and \eqref{pxkkc}, where the form of the transversality condition \eqref{pxkkc} benefits from the ``full" counterpart of the calculus rule in Proposition~\ref{nor_inve} for normals to inverse images that is valid under the full rank assumption in (H3). The measure nonatomicity condition of this theorem is also verified similarly to \cite[Theorem~6.1]{h1}.

Next we establish the new measured coderivative condition \eqref{maximumcondition}, which was not obtained in \cite{h1} even in the particular framework therein. Rewrite first \eqref{gamma_dis} in the form
\begin{eqnarray*}\label{gamma_con}
\gamma^k(t)\in D^*N_{\Th}\big(\psi(\ox^k(\vt^k(t)),\ou^k(\vt^k(t)),\eta^k(t)\big)\Big(\nabla_{x}\psi\big(\ox^k\big(\vt^k(t)\big),\ou^k\big(\vt^k(t)\big)u\Big)\;\mbox{ a.e. }\;t\in[0,T]
\end{eqnarray*}
via the functions $\vt^k(t)$ from \eqref{vt}. Since $\gamma^k(t)$ is a step vector function for each $k\in N$, taking any $t\ne t_j^k$ as $j=0,\ldots,k$ allows us to choose a number $\delta_k>0$ sufficiently small so that $|B|\le\delta_k$ whenever a Borel set $B$ contains $t$, and thus $B$ does not contain any mesh points. Hence $\gamma^k(t)$ remains constant on $B_k$. As a result, we can write the representation
\begin{eqnarray*}
\gamma^k(t)=\frac{1}{|B|}\int_{B}\gamma^k(\t)d\t\;\mbox{ on }\;t\in\B_k\;\mbox{ for all large }\;k\in\N.
\end{eqnarray*}
The separability of the space $C([0,T];\R^s)$ and the boundedness of $\gamma^k(\cdot)$ in $C^*([0,T];\R^s)$ by \eqref{alphak} allow us to select a subsequence of $\{\gg^k(\cdot)\}$ (no relabeling) that weak$^*$ converges in $C^*([0,T];\R^s)$ to some $\gg(\cdot)$. As a result, we get without loss of generality that
\begin{eqnarray*}
\int_B\gamma^k(\t)d\t\to\int_B\gamma(\t)d\t\;\mbox{ as }\;k\to\infty
\end{eqnarray*}
for any Borel set $B$. To proceed further, choose a sequence $\{B_k\}\subset B$ such that $t\in B_k$, $|B_k|\to 0$, and $\frac{1}{|B_k|}\int_{B_k} \gamma(\t)d\t\to\al$ as $k\to\infty$ for some $\alpha\in\R^s$. It follows from the constructions above that
\begin{eqnarray*}
\gamma^k(t)\to\alpha\in\Big(\Limsup_{|B|\to 0}\frac{1}{|B|}\int_B\gamma(\t)d\t\Big)(t)=\Limsup_{|B|\to 0}\frac{\gamma(B)}{|B|}(t).
\end{eqnarray*}
Taking into account the coderivative robustness with respect to all of its variables, we arrive at
\begin{eqnarray*}
\alpha\in D^*N_\Th\big(\psi(\ox(t),\ou(t)),\eta(t)\big)\big(\nabla_x\psi(\ox(t),\ou(t))(q^x(t)-\lm v^x(t))\big)\cap\Limsup_{|B|\to 0}\frac{\gg(B)}{|B|}(t)
\end{eqnarray*}
for a.e.\ $t\in[0,T]$, which verifies the measured coderivative condition \eqref{maximumcondition}.

To complete the proof of the theorem, it remains to justify the nontriviality condition \eqref{nontriv2}. Arguing by contradiction, suppose that \eqref{nontriv2} fails and thus find sequences of $\lambda^k\to 0$ and $\|p_j^k\|\to 0$ as $k\to\infty$ uniformly in $j$. It follows from \eqref{alphak} and \eqref{hamilton2ub} that $\int_0^T\|\gamma^{k}(t)\|dt\to 1$ as $k\to\infty$. Define now the sequence of measurable vector functions $\bb^k\colon[0,T]\to\R^s$ by
\begin{eqnarray*}
\bb^k(t):=\left\{\begin{array}{cl}
\disp\frac{\gamma^k(t)}{\|\gamma^k(t)\|}&{\rm if}\;\gamma^k(t)\ne 0,\\
0&{\rm if}\;\gamma^k(t)=0
\end{array}\right.\quad\mbox{ for all }\;t\in[0,T].
\end{eqnarray*}
Using the Jordan decomposition $\gamma^k=(\gamma^k)^+-(\gamma^k)^-$ gives us a subsequence $\{\gamma^k\}$ and a Borel vector measure $\gamma=\gamma^+-\gamma^-$ such that $\{(\gamma^{k})^+\}$ weak$^*$ converges to $\gamma^+$ and $\{(\gamma^{k})^-\}$ weak$^*$ converges to $\gamma^-$ in $C^*([0,T];\R^s)$. Taking into account the uniform boundedness of $\bb^k(\cdot)$ on $[0,T]$ allows us to apply the convergence result of \cite[Proposition~9.2.1]{vinter} (with $A=A_i:=\B_{\R^s}$ for all $i \in\N$ therein) and thus find a Borel measurable vector functions $\beta^+,\beta^-\colon [0,T]\to\R^s$ so that, up to a subsequence, $\{\beta^k(\gamma^k)^+\}$ weak$^*$ converges to $\beta^+\gamma^+$, and $\{\beta^k(\gamma^k)^-\}$ weak$^*$ converges to $\beta^-\gamma^-$. As a result, we get
\begin{align*}
&\int_0^T\bb^+(t)d\gamma^+(t)-\int_0^T\beta^-(t)d\gamma^-(t)=\lim_{k\to\infty}\int_0^T \beta^k(t)d(\gamma^k)^+(t)-\int_0^T\beta^k(t)d(\gamma^k)^-(t)\\
&=\lim_{k\to\infty}\int_0^T\beta^k(t)d\big((\gamma^k)^+-(\gamma^k)^-\big)(t)=\lim_{k\to\infty}\int_0^T\beta^k(t)d(\gamma^k)(t)=\lim_{k\to\infty}\int_0^T\|
\gamma^{k}(t)\|dt=1.
\end{align*}
This means that $\|\gamma\|=\|\gamma^+\|+\|\gamma^-\|\ne 0$, which contradicts the assumed failure of \eqref{nontriv2}. $\h$\vspace*{0.05in}

Note that the nontrivial optimality conditions obtained in Theorem~\ref{necopt} do not generally exclude the case of their validity for any feasible solution, although even in this case they may be useful as is shown by examples. The following consequence of Theorem~\ref{necopt} presents effective sufficient conditions ensuring {\em nondegenerated} optimality conditions for the considered local minimizers of $(P)$. The reader can find more discussions on nondegeneracy in \cite{arutyunov,AK,vinter} for classical optimal control problems and Lipschitzian differential inclusions and in \cite{h1} for sweepings ones over polyhedral controlled sets.\vspace*{-0.05in}

\begin{Corollary}{\bf(nondegeneracy).}\label{degeneracy} In the setting of Theorem~{\rm\ref{necopt}}, suppose that $\eta(T)$ is well defined and that $\th=0$ is the only vector satisfying the relationships
\begin{eqnarray}\label{int_con_relax}
\th\in D^*N_\Th\big(\psi(\ox(T),\ou(T)),\eta(T)\big)(0),\quad
\nabla\psi\big(\ox(T),\ou(T)\big)^*\th\in\nabla\psi\big(\ox(T),\ou(T)\big)N_\Th\big(\ox(T),\ou(T)\big).
\end{eqnarray}
Then the necessary optimality conditions of Theorem~{\rm\ref{necopt}} hold with the enhanced nontriviality
\begin{eqnarray}\label{enhanced_non}
\lm+{\rm mes}\big\{t\in[0,T]\big|\;q(t)\ne 0\big\}+\|q(0)\|+\|q(T)\|>0.
\end{eqnarray}
\end{Corollary}\vspace*{-0.02in}
{\bf Proof.} Arguing by contradiction, suppose that \eqref{enhanced_non} fails, which yields $\lm=0$, $q(0)=0$, $q(T)=0$ and $q=0$ for a.e.\ $t\in[0,T]$. It follows from \eqref{hamiltonx} that $p\equiv p(T)$ on $[0,T]$. By using $\eqref{q}$ and the fact that $\nabla\psi(\ox(t),\ou(t))$ has full rank on $[0,T]$, we get $\gamma:=\th_1\delta_{\{T\}}+\th_2\delta_{\{T\}}$ for some $\th_1,\th_2\in\R^s$ via the Dirac measures. Let us now check that $\th_1=\th_2=0$. Since $\eta(T)$ is well defined and since $q(T)=0$ and $\lm=0$, we conclude that condition \eqref{maximumcondition} holds at $t=T$ being equivalent to
\begin{eqnarray}\label{condition1}
\th_2\in D^*N_\Th\big(\psi(\ox(T),\ou(T)),\eta(T)\big)(0).
\end{eqnarray}
On the other hand, it follows from $q(T)=0$ due to \eqref{q} that $p(T)=\nabla\psi(\ox(t),\ou(t))^*\th$. Using further \eqref{pxkkc} tells us that $\nabla\psi(\ox(T),\ou(T))^*\th_2\in-\nabla\phi(\ox(T),\ou(T))N_\Th((\ox(T),\ou(T))$. Then it follows from \eqref{int_con_relax} and
\eqref{condition1} that $\th_2=0$, which yields $\th_1=0$ and gives us a contradiction. $\h$\vspace*{-0.07in}

\begin{Remark}{\bf(discussions on nondegeneracy).}\label{example_nondegenerate} {\rm It is easy to see that the imposed assumption \eqref{enhanced_non} excludes the degeneracy case of $\lm:=0$, $p=q\equiv 0$, and $\gamma:=\delta_{\{T\}}$ in Theorem~\ref{necopt}. Furthermore, the inferiority assumption $\psi(\ox(T),\ou(T))\in{\rm int}\,\Th$ yields \eqref{int_con_relax} while not vice versa. To illustrate it, consider the following {\em example}: minimize the cost
\begin{eqnarray*}
J[u]:=\frac{1}{2}\Big(x(2)-1\Big)^2+\int_0^1\Big(u(t)+2-t\Big)^2dt+\int_1^2\Big(u(t)+1\Big)^2dt
\end{eqnarray*}
over the dynamics $-\dot x(t)\in N(x(t);(-\infty,u])$ with $x(0)=3/2$ and $u(0)=-2$. We can directly check that the only optimal trajectory in this problem is given by $\ox(t)=3/2$ on $[0,1/2]$, $\ox(t)=2-t$ on $[1/2,1]$, and $\ox(t)=1$ on $[1,2]$, It is generated by the optimal control $\ou(t)=t-2$ on $[0,1]$ and $\ou(t)=-1$ on $(1,2]$.

To check the nondegeneracy condition \eqref{int_con_relax} in this example with $\psi(x,u)=x+u$ and $\Th=(-\infty,0]$, observe that the second inclusion therein reduces to
\begin{eqnarray*}
(\al,\al)\in-N\big((1,-1);\{(x,u)|\,x+u\le 0\}\big),
\end{eqnarray*}
which is equivalent to $\al\le 0$. The first inclusion reads as
\begin{eqnarray*}
\al\in D^*N_{\R_-}\big(0,0\big)\big(0\big)=[0,\infty)
\end{eqnarray*}
giving us $\al\ge 0$. Thus $\al=0$, and condition \eqref{int_con_relax} is satisfied. On the other hand, we see that the point $(\ox(2),\ou(2))=(1,-1)$ does not belong to the interior of the set $\Th$.}
\end{Remark}\vspace*{-0.25in}

\section{Hamiltonian Formalism and Maximum Principle}\label{hamiltonian}
\setcounter{equation}{0}

The necessary optimality conditions for $(P)$ obtained in Theorem~\ref{necopt} and Corollary~\ref{degeneracy} are of the extended {\em Euler-Lagrange} type that is pivoting in optimal control of Lipschitzian differential inclusions; see \cite{mordukhovich,vinter}. The result of \cite[Theorem~1.34]{mordukhovich} tells us that the Euler-Lagrange framework involving coderivatives implies the {\em maximum condition} of the Weierstrass-Pontryagin type for problems with convex velocities provided that the velocity mapping is {\em inner semicontinuous} (e.g., Lipschitzian), which is never the case in our setting \eqref{F}. Nevertheless, we show in what follows that the Hamiltonian formalism and the maximum condition can be derived from the measured coderivative condition of Theorem~\ref{necopt} in rather broad and important situations by using {\em coderivative calculations} available in variational analysis.

The first result of this section deals with problem $(P)$ in the case where $\Th:=\R^s_-$. In this case we consider the set of {\em active constraint indices}
\begin{eqnarray}\label{active_constr}
I(x,u):=\big\{i\in\{1,\ldots,s\}\big|\;\psi_i(x,u)=0\big\}.
\end{eqnarray}
It follows from Proposition~\ref{nor_inve} under (H3) that for each $v\in-N(x;C(u))$ there is a unique collection $\{\alpha_i\}_{i\in I(x,u)}$ with $\alpha_i\le 0$ and $v=\sum_{i\in I(x,u)}\alpha_i[\nabla_x\psi(x,u)]_i$. Given $\nu\in\R^s$, define the vector $[\nu,v]\in\R^n$ by
\begin{eqnarray}\label{nu}
[\nu,v]:=\sum_{i\in I(x,u)}\nu_i\alpha_i\big[\nabla_x\psi(x,u)\big]_i
\end{eqnarray}
and introduce the {\em modified Hamiltonian} function
\begin{eqnarray}\label{new_hamiltonian}
H_{\nu}(x,u,p):=\sup\big\{\big\la[\nu,v],p\big\ra\big|\;v\in-N\big(x;C(u)\big)\big\},\quad(x,u,p)\in\R^n\times\R^m\times\R^n.
\end{eqnarray}

The following consequence of Theorem~\ref{necopt} and Corollary~\ref{degeneracy} shows that the measured coderivative condition \eqref{maximumcondition} yields the {\em maximization} of the modified Hamiltonian \eqref{new_hamiltonian} at the local optimal solutions for $(P)$ with polyhedral constraints corresponding to $\Th=\R^s_-$ in \eqref{mov-set}.\vspace*{-0.05in}

\begin{Corollary}{\bf(maximum condition for polyhedral sweeping control systems under surjectivity).}\label{maximumconditionorthant} In the frameworks of Theorem~{\rm\ref{necopt}} and Corollary~{\rm\ref{degeneracy}} with $\Th=\R^s_-$ we have the corresponding necessary optimality conditions therein together with the following maximum condition: there is a measurable vector function $\nu\colon[0,T]\to\R^s$ such that $\nu(t)\in\Limsup_{|B|\to 0}\frac{\gamma(B)}{|B|}(t)$ and
\begin{eqnarray}\label{novelmaxcond}
\big\langle\big[\nu(t),\dot\ox(t)\big],q^x(t)-\lm v^x(t)\big\rangle=H_{\nu(t)}\big(\ox(t),\ou(t),q^x(t)-\lm v^x(t)\big)=0\;\mbox{ a.e. }\;t\in[0,T].
\end{eqnarray}
\end{Corollary}\vspace*{-0.05in}
{\bf Proof.} Let us show that the maximum condition \eqref{novelmaxcond} follows from the measured coderivative condition \eqref{maximumcondition}. To proceed, we need to compute the coderivative of the normal cone mapping $D^*N_{\R^s_-}$, which has been done in variational analysis in several forms; see, e.g., \cite{mord,BR} for more discussions and references. We use here the one taken from \cite{MO07}: if $D^*N_{\R^s_-}(w,\xi)(u)\ne\emp$, then
\begin{eqnarray}\label{orthant}
D^*N_{\R^s_-}(w,\xi)(u)=\left\{\omega\in\mathbb{R}^{s}\left\vert\begin{array}{ll}
\omega_i=0&\mbox{if }\;w_i<0\;\mbox{ or if }\;w_i=0,\;\xi_{i}=0,\;u_i<0
\\\omega_i\ge 0&\mbox{if }\;w_i=0,\;\xi_{i}=0,\;u_i\ge 0\\
\omega_i\in\R&\mbox{if }\;\xi_{i}>0,\;u_i=0
\end{array}\right.\right\}.
\end{eqnarray}
The measured coderivative condition \eqref{maximumcondition} reads in this case as:
\begin{eqnarray*}
\nu(t)\in D^*N_{\R^s_-}\big(\psi(\ox(t),\ou(t)),\eta(t)\big)\big(\nabla_x \psi(\ox(t),\ou(t))(q^x(t)-\lm v^x(t))\big)\;\mbox{ a.e. }\;t\in[0,T]
\end{eqnarray*}
with a vector function $\nu(t)\in\Limsup_{|B|\to 0}\frac{\gamma(B)}{|B|}(t)$, which can be selected as (Lebesgue) measurable on $[0,T]$ due to the well-known measurable selection results; see, e.g., \cite{rw,vinter}. It follows from \eqref{orthant} that
\begin{eqnarray}\label{eta}
\big[\eta_{i}(t)>0\big]\Longrightarrow\Big[\big\la\lm v^x(t)-q^{x}(t),\big[\nabla_x\psi\big(\ox(t),\ou(t)\big)\big]_i\big\ra=0\Big]\;\mbox{ for a.e. }\;t\in[0,T],\quad i=1,\ldots,s,
\end{eqnarray}
which gives us by equation \eqref{etajkl} that
\begin{eqnarray}\label{LHS-MX}
\big\la\big[\nu(t),\dot\ox(t)\big],\lm v^x(t)-q^{x}(t)\big\ra=0\;\mbox{ a.e. }\;t\in[0,T].
\end{eqnarray}
On the other hand, we get from \eqref{active_constr}--\eqref{new_hamiltonian} with $I:=I(\ox(t),\ou(t))$ that
\begin{eqnarray}\label{RHS-MX}
H_{\nu(t)}\big(\ox(t),\ou(t),q^x(t)-\lm v^x(t)\big)=\sup\Big\{\sum_{i\in I}\alpha_i\nu_i(t)\big\langle\big[\nabla_x\psi\big(\ox(t),\ou(t)\big)\big]_i,q^x(t)-\lm v^x(t)\big\rangle|\;\alpha_i\le 0\Big\}.
\end{eqnarray}
for a.e.\ $t\in[0,T]$. Applying now \eqref{orthant} gives us the implication
\begin{eqnarray*}
\big[i\in I\big(\ox(t),\ou(t)\big)]\Longrightarrow\Big[\nu_i(t)\big\langle\big[\nabla_x\psi\big(\ox(t),\ou(t)\big)\big]_i,q^x(t)-\lm v^x(t)\big\rangle\ge 0\Big]\;\mbox{ a.e. }\;t\in[0,T].
\end{eqnarray*}
Combining this with \eqref{RHS-MX}, we get $H_{\nu(t)}(\ox(t),\ou(t),q^x(t)-\lm v^x(t))=0$ for a.e.\ $t\in[0,T]$ and thus arrive at the maximum condition \eqref{novelmaxcond}, where the other equality was established in \eqref{LHS-MX}. $\h$\vspace*{0.05in}

Observe that the explicit coderivative computation \eqref{orthant} plays a crucial role in deriving the maximum condition \eqref{novelmaxcond} in Corollary~\ref{maximumconditionorthant}. Available second-order calculus and coderivative evaluations for the normal cone mappings allow us to derive more
general results of the maximum principle type in sweeping optimal control. The next theorem addresses the case where the set $\Th$ in \eqref{mov-set} is given by
\begin{eqnarray}\label{h}
\Th=h^{-1}(\R^l_-):=\big\{z\in\R^s\big|\;h(z)\in\R^l_-\big\}
\end{eqnarray}
via a smooth mapping $h\colon\R^s\to\R^l$. As mentioned above, the surjectivity condition on the Jacobian $\nabla h(\oz)$ at a fixed point $\oz$ corresponds to the LICQ condition at $\oz$. Dealing with linear mappings $h(z):=Az-b$, we may replace the LICQ condition in the coderivative evaluation by a weaker {\em positive LICQ} (PLICQ) condition at $\ox$ that is discussed and implemented in \cite{h1}. It is used in what follows.\vspace*{-0.07in}

\begin{Theorem}{\bf(maximum principle in sweeping optimal control).}\label{max-pr} Consider the control problem $(P)$ in the frameworks of Theorem~{\rm\ref{necopt}} and Corollary~{\rm\ref{degeneracy}} with the set $\Th$ given by \eqref{h}, where $h\colon\R^s\to\R^l$ is ${\cal C}^2$-smooth around the local optimal solution $\oz(t):=(\ox(t),\ou(t))$ for all $t\in[0,T]$. Suppose that either $\nabla h(\oz(t))$ is surjective, or $h(\cdot)$ is linear and the PLICQ assumption is fulfilled at $\oz(t)$ on $[0,T]$. Then, in addition to the corresponding necessary optimality conditions of the statements above, the maximum condition \eqref{novelmaxcond} holds with a measurable vector function $\nu\colon[0,T]\to\R^s$ satisfying the inclusion
\begin{eqnarray}\label{tildegg1}
\nu(t)\in D^*N_{\R^l_-}\big(h(\psi(\ox(t),\ou(t))),\mu(t)\big)\big(\nabla_x\psi(\ox(t),\ou(t))(q^x(t)-\lm v^x(t))\big)\;\mbox{ a.e. }\;t\in[0,T],
\end{eqnarray}
where $\mu\colon[0,T]\to\R^l$ is also measurable and such that
\begin{eqnarray*}\label{v}
\mu(t)\in N_{\R^l_-}\big(h(\psi(\ox(t),\ou(t))\big)\;\mbox{ with }\;\eta(t)=\nabla h\big(\psi(\ox(t),\ou(t)\big)^*\mu(t)\;\mbox{ a.e. }\;t\in[0,T].
\end{eqnarray*}
\end{Theorem}\vspace*{-0.05in}
{\bf Proof.} As in the proof of Corollary~\ref{maximumconditionorthant}, we derive from the measured coderivative condition \eqref{maximumcondition} the existence of a measurable function $\Tilde\nu\colon[0,T]\to\R^l$ satisfying the inclusion
\begin{eqnarray}\label{v1}
\Tilde\nu(t)\in D^*N_{h^{-1}(\R^l_-)}\big(\psi(\ox(t),\ou(t)),\eta(t)\big)\big(\nabla_x\psi(\ox(t),\ou(t))(q^x(t)-\lm v^x(t))\big)\;\mbox{ a.e. }\;t\in[0,T].
\end{eqnarray}
Assuming that $\nabla h(\oz(t))$ is surjective on $[0,T]$ and applying the second-order chain rule from \cite[Theorem~1.127]{mordukhovich} together with the aforementioned measurable selection results, we find measurable functions $\nu\colon[0,T]\to\R^s$ and $\mu\colon[0,T]\to\R^l$ satisfying the conditions \eqref{tildegg1} and \eqref{v1} as well as
\begin{eqnarray*}
\Tilde\nu(t)=\nabla^2\big\langle\mu(t),h\big\rangle\big(\psi(\ox(t),\ou(t))\big)+\nabla h\big(\psi(\ox(t),\ou(t))\big)^*\nu(t),
\end{eqnarray*}
which uniquely determines $\nu(t)$ from $\Tilde\nu(t)$ for a.e.\ $t\in[0,T]$. The validity of the maximum condition follows now from the proof of Corollary~\ref{maximumconditionorthant} due to the relationships above.

In the case where $h$ is linear and the PLICQ holds, we proceed similarly by applying the evaluation of $D^*N_{h^{-1}(\R^l_-)}$ from \cite[Lemma~4.2]{h1} without claiming that $\nu(t)$ is uniquely defined by $\Tilde\nu(t)$. $\h$\vspace*{-0.05in}

\begin{Remark}{\bf(discussions on the maximum principle).}\label{max-disc} {\rm The necessary optimality conditions of the maximum principle type obtained in Corollary~\ref{maximumconditionorthant} and Theorem~\ref{max-pr} are the first in the literature for sweeping processes with controlled moving sets. Note that our form of the modified Hamiltonian \eqref{new_hamiltonian} is different from the {\em conventional Hamiltonian} form
\begin{eqnarray}\label{hamilton}
H(x,p):=\sup\big\{\la p,v\ra\ra\big|\;v\in F(x)\big\}
\end{eqnarray}
used in optimal control of Lipschitzian differential inclusions $\dot x\in F(x)$ that extends the Hamiltonian in classical optimal control. We show below in Example~\ref{counterexample} that the maximum principle via the conventional Hamiltonian \eqref{hamiltonian} {\em fails} for our problem $(P)$. The reason is that \eqref{hamilton} does not reflects {\em implicit} state constraints, which do not appear for Lipschitzian problems while being an essential part of the sweeping dynamics; see the discussions in Section~1. Note that the Hamiltonian form \eqref{new_hamiltonian} is also different from those used in \cite{ac,bk,pfs} for deriving maximum principles in sweeping control problems with uncontrolled moving sets that are significantly diverse from our problem $(P)$.

We can see that the new maximum principle form \eqref{novelmaxcond} incorporates vector measures that appear through the measured coderivative condition \eqref{maximumcondition}. The fact that measures naturally arise in descriptions of necessary optimality conditions in optimal control problems with state constraints has been first realized by Dubovitskii and Milyutin \cite{duboviskii} and since that has been fully accepted in the literature; see, e.g., \cite{arutyunov,vinter} and the references therein. There are interesting connections between our form of the maximum principle for controlled sweeping processes and the Hamiltonian formalism in models of contact and nonsmooth mechanics (see, e.g., \cite{brog,razavy}), which we plan to fully investigate in subsequent publications.}
\end{Remark}\vspace*{-0.07in}

The next example shows that the maximum principle in the conventional form used for Lipschitzian differential inclusions \cite{arutyunov,vinter} with the standard Hamiltonian \eqref{hamilton} fails for the sweeping control problem $(P)$, while our new form obtained in Theorem~\ref{max-pr} holds.\vspace*{-0.05in} 

\begin{Example}{\bf(failure of the conventional maximum principle for sweeping processes).}\label{counterexample} {\rm We consider the optimal control problem for the sweeping process taken from \cite[Example~7.6]{h1}, where the controlled moving set is defined by \eqref{sw-con2} with control actions $u_j(t)$ and $b_j(t)$. Specify the initial data as
\begin{eqnarray*}
n=m=2,\;x_0=(1,1),\;T=1,\;\varphi(x)=\frac{\|x\|^2}{2},\;\mbox{ and }\;\ell(t,x,u,b,\dot x,\dot u,\dot b):=\frac{1}{2}\big(\dot{b}_1^2+\dot{b}_2^2\big)
\end{eqnarray*}
and fix the $u$-controls as $\ou_1\equiv(1,0)$, $\ou_2\equiv(0,1)$. The necessary optimality conditions of Corollary~\ref{degeneracy} give us in this case the following relationships on $[0,1]$:
\begin{itemize}
\item[(1)] $w(\cdot)=0$, $v^x(\cdot)=0$, $v^b(\cdot)=\big(\dot{b}_1(\cdot),\dot{b}_2(\cdot)\big)$;\quad(2)$\;\dot{\ox}_i(t)\ne 0\Longrightarrow
q^x_i(t)=0$,\;$i=1,2$,
\item[(3)] $p^b(\cdot)$ is constant with nonnegative components; $-p_i^x(\cdot)=\lambda\ox(1)+p_i^b(\cdot)\ou_i$ are constant for $i=1,2$,
\item[(4)] $q^x(t)=p^x-\gamma([t,1])$,\quad $q^b(t)=\lambda\dot{\ob}(t)=p^b+\gamma([t,1])$ for a.e.\ $t\in[0,1]$,
\item[(5)] $\lambda+\|q(0)\|+\|p(1)\|\ne 0$ with $\lm\ge 0$.
\end{itemize}
Observe first that the pair $\ox(t)=(1,1)$ and $\ob(t)=0$ on $[0,1]$ satisfies the necessary conditions with $p_1^x=p_2^x=-1$, $p_1^b=p_2^b=\gamma_1=\gamma_2=0$, and $\lambda=1$). The conventional Hamiltonian \eqref{hamilton} reads now as
\begin{eqnarray*}
H(x,b,p)=\sup\big\{\la p,v\ra\big|\;v\in-N\big(x;C((1,0),(0,1)),b\big)\big\},
\end{eqnarray*}
and we get by the direct calculation that
\begin{align*}
H\big(\ox(t),\ob(t),q^x(t)-\lambda v^x(t)\big)=&H\big((1,1),(1,1),(-1,-1)\big)\\
=&\sup\big\{\big\la(-1,-1),v\big\ra\big|\;v\in-N\big((1,1);C(((1,0),(0,1)),(1,1)\big)\big\}\\
=&\sup\big\{\big\la(-1,-1),v\big\ra\big|\;v_1\le 0,\;v_2\le 0\big\}=\infty 
\end{align*}
while $\la\dot\ox(t),q^x(t)-\lambda v^x(t)\ra=0$, and thus the conventional maximum principle fails in this example. At the same time, the new maximum condition \eqref{novelmaxcond} holds trivially with $\nu(t)\equiv 0$ on $[0,1]$.}
\end{Example}\vspace*{-0.05in}

The following consequence of Theorem~\ref{max-pr} provides a natural sufficient condition for the validity of the maximum principle in terms of the conventional Hamiltonian \eqref{hamilton}.\vspace*{-0.05in}

\begin{Corollary}{\bf(maximum principle in the conventional form).}\label{con-mp} Assume that in the setting of Theorem~{\rm\ref{max-pr}} we have the condition
\begin{eqnarray*}
\eta_i(t)>0\;\mbox{ for all }\;i\in I\big(\ox(t),\ou(t)\big)\;\mbox{ a.e. }\;t\in[0,T].
\end{eqnarray*}
where $\eta(t)\in N(\psi(\ox(t),\ou(t));\Th)$ is uniquely defined by \eqref{etajkl}. Then
\begin{eqnarray*}
\big\langle\dot\ox(t),q^x(t)-\lm v^x(t)\big\rangle=H\big(\ox(t),\ou(t),\ob(t),q^x(t)-\lm v^x(t)\big)=0\;\mbox{ a.e. }\;t\in[0,T].
\end{eqnarray*}
\end{Corollary}\vspace*{-0.05in}
{\bf Proof.} This follows from \eqref{eta} and its counterpart in the proof of Theorem~\ref{max-pr} by definitions of the modified and conventional Hamiltonians. $\h$\vspace{-0.15in}

\section{Applications to Elastoplasticity and Hysteresis}\label{application-example}
\setcounter{equation}{0}

In this section we discuss some applications of the obtained necessary optimality conditions to a fairly general class of problems relating to elastoplasticity hysteresis. 

Let us consider the model of this type discussed in \cite[Section~3.2]{adly}, which can be described in our form, where $Z$ is a closed convex subset of the $\frac{1}{2}n(n+1)$-dimensional vector space $E$ of symmetric tensors $n\times n$ with ${\rm int}\,Z\ne\emp$. Using the notation of \cite{adly}, define the strain tensor $\epsilon=\{\epsilon\}_{i,j}$ by $\epsilon:=\epsilon^e+\epsilon^p$, where $\epsilon^e$ is the elastic strain and $\epsilon^p$ is the plastic strain. The elastic strain $\epsilon^e$ depends on the stress tensor $\sigma=\{\sigma\}_{i,j}$ linearly, i.e., $\epsilon^e=A^2\sigma$, where $A$ is a constant symmetric positive-definite matrix. The {\em principle of maximal dissipation} says that
\begin{eqnarray}\label{elasto-plastic}
\big\langle\dot\epsilon^p(t),z-\sigma(t)\big\rangle \le 0 \;\mbox{ for all }\;z\in Z.
\end{eqnarray}
It is shown in \cite{adly} that the variational inequality \eqref{elasto-plastic} is equivalent to the {\em sweeping processes}
\begin{eqnarray}\label{sweeping-elasto}
\dot\zeta(t)\in-N\big(\zeta(t);C(t)\big),\;\zeta(0)=A\sigma(0)-A^{-1}\epsilon(0)\in C(0),
\end{eqnarray}
where $\zeta(t):=A\sigma(t)-A^{-1}\epsilon(t)$ and $C(t):=-A^{-1}\epsilon(t)+AZ$. It can be rewritten in the frame of our problem $(P)$ with $x:=\zeta$, $u:=\epsilon$, $\psi(x,u):=x+A^{-1}u$, and $\Th:=AZ$. Thus we can apply Theorem~\ref{necopt} to this class of hysteresis operators for the {\em general} elasticity domain $Z$. Note that a similar model is considered in \cite{herzog1}, only for the von Mises yield criterion. Our results obtained here give us the flexibility of applications to many different elastoplasticity models including those with the Drucker-Prager, Mohr-Coulomb, Tresca, von Mises yield criteria, etc.

In the following example we summarize applications of Theorem~\ref{necopt} to solve a meaningful optimal control problem generated by the elastoplasticity dynamics \eqref{sweeping-elasto}.\vspace*{-0.05in}

\begin{Example}{\bf(optimal control in elastoplasticity).}\label{example1} {\rm Consider the dynamic optimization problem:
\begin{eqnarray*}
\mbox{minimize }\;J[\epsilon]:=\int_0^1\frac{1}{2}\|\dot\epsilon(t)\|^2dt+\frac{1}{2}\|\zeta(1)-\zeta_1\|^2
\end{eqnarray*} 
over feasible solutions to the sweeping process \eqref{sweeping-elasto} with the initial point $\zeta(0)\in C(0)$. Observe that the linear function $\bar\epsilon(t):=tz+\bar\epsilon(0)$ with appropriate adjustments of the starting and terminal points is an optimal control to the corresponding problem $(P)$. Remembering the above notation for $x,u,\psi$, and $\Th$, we derive from Theorem~\ref{necopt} the following necessary optimality conditions:
\begin{eqnarray*}
\begin{array}{ll}
(1)\;\;\dot p(t)=\lm (0,0),\quad(2)\;\;q^u(t)=\lm \frac{d}{dt}\bar \epsilon(t)\;\mbox{ a.e. }\;t\in[0,T],\\
(3)\;\;q(t)=p(t)-\disp\int_{[t,T]}(1,1)d\gamma(s)\;\mbox{ a.e. }\;t\in[0,T],\\
(4)\;\;-\big(p^x(1),p^b(1)\big)=\lambda\big(\big(\zeta(1)-\zeta_1\big),0\big)+N\big((\zeta(1),\epsilon(1));\{(\zeta,\epsilon)|\,\zeta+A^{-1}\epsilon\in AZ\}\big),\\
(5)\;\;\lm +\sup_{t\in[0,T]}\|p(t)\|+\|\gamma\|\ne 0.
\end{array}
\end{eqnarray*}
Assuming that $\zeta(T)\in{\rm int}\,C(T)$ and using the measure nonatomicity condition, we get from (1)--(5) that $\lm=1$. It follows from the linearity of $\bar\epsilon$ that $\frac{d}{dt}\bar\epsilon(t)\equiv:a$ and that the choice of $p\equiv 0$, $q\equiv(\lm a,\lm a)$, and $\gamma \equiv(-a,-a)\delta_{\{1\}}$ fulfills all the conditions (1)--(5).}$\h$
\end{Example}\vspace*{-0.05in}

Note that in the above example we actually guessed the form of optimal solutions and then checked the fulfillment of the obtained necessary optimality conditions. The following example does more by showing how to calculate an optimal solution by using the necessary optimality conditions of Theorem~\ref{necopt} together with the maximum condition from Corollary~\ref{maximumconditionorthant}.\vspace*{-0.05in}

\begin{Example}{\bf(calculation of optimal solutions).}\label{Ex2} {\rm Consider here the optimal control problem taken from the example in Remark \ref{example_nondegenerate}. Applying the necessary optimality conditions of Corollary~\ref{degeneracy} with the enhanced nontriviality condition \eqref{enhanced_non}, we get the following:
\begin{eqnarray*}
\begin{array}{ll}
(1)\;\;w=0,\;v(t)=\big(0,2(\ou(t)+2-t)\big)\;\mbox{ if }\;t\in[0,1]\;\mbox{ and }\;v(t)=\big(0,2(\ou(t)+1)\big)\;\mbox{ if }\;t\in(1,2],\\
(2)\;\;\dot{\ox}(t)\in-N\big(\ox(t);(-\infty,-\ou(t)]\big)\;\mbox{ a.e. }\;t\in[0,2],\\
(3)\;\big(\dot{p}^x,\dot{p}^u\big)(t)=\big(0,0\big)\;\mbox{ a.e. }\;t\in[0,2],\\
(4)\;\;q^u(t)=2\lm\big(\ou(t)+2-t\big)\;\mbox{ if }\;t\in[0,1]\;\mbox{ and }\;q^u(t)=2\lm\big(\ou(t)+1\big)\;\mbox{ if }\;t\in(1,2],\\
(5)\;(q^x,q^u)(t)=(p^x,p^u)(t)-\disp\Big(\int_t^2d\gamma,\int_t^2d\gamma\Big)\;\mbox{ a.e. }\;t\in[0,2],\\
(6)-\big(p^x(2),p^u(2)\big)=\lambda\big(\big(\ox(2)-1\big),0\big)+N\big((\ox(2),\ou(2));\{(x,u)|\,x+u\le 0\}\big),\\
(7)\;\;\lm +{\rm mes}\big\{t\in[0,2]\big|\;q(t)\ne 0\big\}+\|q(0)\|+\|q(T)\|>0,
\end{array}
\end{eqnarray*}
together with the measured coderivative condition telling is that
\begin{eqnarray}\label{MXC}
\nu(t)\in D^*N_{\R_-}\big(\ox(t)+\ou(t),-\dot \ox(t)\big)\big(q^x(t)\big)\;\mbox{ a.e. }\;t\in[0,2],
\end{eqnarray}
where $\nu(t)\in\Limsup_{|B|\to 0}\frac{\gamma(B)}{|B|}(t)$. To proceed, take $\sigma\in(0,2]$ such that the state constraint is inactive on $[0,\sigma)$ while it is active on $[\sigma,2]$. It follows from (3) that $(p^x,p^u)$ is constant on $[0,2]$. Assuming by contradiction that $\lm=0$, we get from (4) that $q^u\equiv 0$, and so $q^x\equiv 0$ by (5). It shows that (7) is violated, and thus $\lm=1$. Since $\ox(t)+\ou(t)<0$ on $[0,\sigma)$, it follows from \eqref{MXC} and \eqref{orthant}, which gives us the maximum condition \eqref{novelmaxcond}, that $\nu(t)=0$ for all $t\in[0,\sigma)$. Then $q^u(t)$ remains constant on $[0,\sigma)$ by (5). Combining this with (4) shows that $\ou(t)=t-2$ on $[0,\sigma]$. Using now (2) and the fact that
$\ox(t)+\ou(t)<0$ for $t\in[0,\sigma)$, we get $\ox\equiv 3/2$ on $[0,\sigma]$. Then $\ox(\sigma)+\ou(\sigma)=\sigma-2+3/2=0$ and $\sigma=1/2$. Assuming further that $\dot\ox(t)<0$ on $(1/2,\sigma_1)$ for some $\sigma_1\in(1/2,2]$ and using \eqref{MXC} together with \eqref{orthant} tell us that $q^x\equiv 0$ on $(1/2,\sigma_1)$. Applying (5) again, we have $\ou(t)=t-2$ on $[1/2,\sigma_1]$. Since $\dot\ox(t)<0$ on $(1/2,\sigma_1)$, it follows from (2) that $\ox(t)=-\ou(t)=2-t$ on $[1/2,\sigma_1]$; thus we tend to move $\ox(t)$ towards $1$. This means that when $\ox(\sigma_1)=2-\sigma_1=1$, the object stops moving, i.e., $\dot\ox(t)=0$ on $[1,2]$. Note that (5) yields $\nu(t)=-\dot q^u(t)=-\dot q^x(t)$ on $[0,2]$. If furthermore $\ou(t)$ is strictly increasing, then $q^u(t)$ is also strictly increasing. Since $\nu(t)=-\dot q^u(t)$, we get $\nu(t)<0$. Then \eqref{MXC} and \eqref{orthant} imply that $\dot\ox(t)<0$, a contradiction. Assuming that $\ou(t)$ is strictly decreasing tells us by (4) that $q^u(t)$ is strictly decreasing, and so $\nu(t)<0$ by $\nu(t)= -\dot q^u(t)$. In this case conditions \eqref{MXC} and \eqref{orthant} yield $q^x(t)\ge 0$, and hence we have $q^x(2)\le 0$ by (6). It shows that $q^x(t)=0$ on $[1,2]$ and so $\nu(t)=-\dot q^x(t)=0$, which contradicts the condition $\nu(t)<0$. Thus $\ou(t)$ remains constant on $[1,2]$, and we find an optimal solution to  this problem.}
\end{Example}\vspace*{-0.25in}

\section{Conclusions}\label{conclusion}
\setcounter{equation}{0}

This paper presents new results on extended Euler-Lagrange and Hamiltonian optimality conditions for a rather general class of controlled sweeping processes. In our future research in this direction we plan to focus on optimal control problems governed by {\em rate-independent operators} having the following description. Given two functionals $E\colon[0,T]\times Z\to\R$ and $\Psi\colon Z\times Z\to[0,\infty)$ on a Banach (or finite-dimensional) space $Z$, we consider the {\em doubly nonlinear evolution inclusion}
\begin{eqnarray*}
0\in\partial_v\Psi\big(z(t),\dot z(t)\big)+\partial E\big(t,z(t)\big)\;\mbox{ a.e. }\;t\in[0,T].
\end{eqnarray*}
If $E$ is smooth, the inclusion above is equivalent to
\begin{eqnarray*}
\dot z(t)\in N_{C(z(t))}\big(\nabla E(t,z(t))\big)\mbox{ a.e. }\;t\in[0,T],
\end{eqnarray*}
where $\{C(z)\}_{z\in Z}$ is the family of closed convex subsets of $Z$ related to $\Psi$ by the formula
\begin{eqnarray*}
\Psi(z,v):=\sup\big\{(\sigma,v)\big|\;\sigma\in C(z)\big\}\;\mbox{ for all }\;z,v\in Z.
\end{eqnarray*}
Among our major applications we plan to consider practical hysteresis models, especially those arising in problems of contact and nonsmooth mechanics; see \cite{brog,razavy}.
\end{document}